\documentclass[12pt]{article}
\usepackage[margin=1in]{geometry}
\usepackage{graphicx}          
\usepackage{bbm}
\usepackage[normalem]{ulem}
\usepackage[usenames]{color}

\usepackage{times}
\usepackage{textcomp}
\usepackage{cite}
\usepackage{url}
\usepackage{subcaption}
\usepackage{amsfonts,mathrsfs}
\usepackage{amssymb,amsmath}
\usepackage{hyperref}
\usepackage{amsthm}

\usepackage{verbatim}
\usepackage{acronym}
\usepackage{mathtools}

\usepackage{todonotes}

\usepackage{enumitem}
\usepackage{subfiles}
\usepackage{algorithm}
\usepackage{algpseudocode}
\usepackage{hyperref}

\newtheorem{thm}{Theorem}

\newtheorem{assum}{Assumption}

\newtheorem{conj}[thm]{Conjecture}
\newtheorem{cor}{Corollary}

\newtheorem{defn}{Definition}

\newtheorem{lem}{Lemma}

\newtheorem{prob}[thm]{Problem}

\newenvironment{pf}%
  {\par\addvspace{\baselineskip}\noindent
   {\bfseries \Elproofname}\enspace\ignorespaces}%
  {\hfill$\qedsymbol$\par\addvspace{\baselineskip}}
\def\Elproofname{Proof.}
\newenvironment{pf*}[1]%
  {\par\begingroup\def\Elproofname{#1}%
   \pf}%
  {\endpf\endgroup}

\author{Ashwin Verma, Soheil Mohajer, and Behrouz Touri\footnote{
A.\ Verma and B.\ Touri are with the ECE Department of University of California San Diego (email: \{a1verma, btouri\}@ucsd.edu), and S.\ Mohajer is with the ECE Department of the University of Minnesota (email: \{soheil\}@umn.edu). 
This research is supported by AFOSR under the grant FA9550-23-1-0057. We thank Alireza Sharbafchi for his helpful comments and discussions while completing this work. }}
\title{Multi-Agent Fact Checking} 
\date{}

\definecolor{OliveGreen}{rgb}{0,0.5,0}
\newcommand{\av}[1]{{\color{black} #1}}
\newcommand{\avnew}[1]{{\color{black} #1}}
\newcommand{\beh}[1]{{\color{black} #1}}

\newcommand{\soh}[1]{{\color{black} #1}}

\usepackage{tikz}
\usetikzlibrary{patterns}
\newcommand{\lpi}[1]{\ell_{\pi_{#1}}}
\definecolor{ao}{rgb}{0.0, 0.5, 0.0}

 \usepackage{bbm}

\usepackage{ulem}

\def \E {\mathbb{E}}
\def \N {\mathbb{N}}

\def \R {\mathbb{R}}

\def \Ecal {\mathcal{E}}
\def \Fcal {\mathcal{F}}

\def \Ncal {\mathcal{N}}

\def \Scal {\mathcal{S}}
\def \Wcal {\mathcal{W}}
\def \Wbar {\bar{\Wcal}}

\def \one {\vec{1}}

\def \Scal {\mathcal{S}}

\def \bpi {\Bar{\pi}}

\newcommand {\indicator}[1]{\mathbbm{1}_{\{#1\}}}

\renewcommand{\vec}[1]{\boldsymbol{#1}}

\def \Pr{\mathbb{P}}
\newcommand\inp[2]{\left\langle #1, #2 \right\rangle}

\usepackage{ulem}

\newcommand{\vx}{\vec{x}}

\newcommand{\vr}{\vec{r}}

\newcommand{\ftilde}{{\tilde f}}
\newcommand{\vftilde}{\vec{\tilde f}}
\newcommand{\Kcal}{\mathcal{K}}
\newcommand{\proj}[1]{\gamma(#1)}

\newcommand{\paran}[1]{\left(#1\right)}
\newcommand{\norm}[1]{\left\| #1 \right\|}

\newcommand{\f}{\tilde f}

\newcommand{\transpose}[1]{{#1}^\top}

\newcommand{\bl}{\boldsymbol{\ell}}

\newcommand{\vf}{\vec{f}}
\newcommand{\vP}{\vec{P}}
\newcommand{\vPpr}{\vec{P}_{\text{pr}}}
\newcommand{\hvPpr}{\hat{\vec{P}}_{\text{pr}}}

\newcommand{\vM}{\vec{M}}
\newcommand{\vR}{\vec{R}}

\newcommand{\vpi}{\vec{\pi}}

\newcommand{\stepq}{\nu}

\newcommand{\vcorrection}[1]{\vec{\rho}({#1})}

\newcommand{\Dkl}{D_{\text{KL}}}
\newcommand{\del}[2]{\frac{\partial #1}{\partial #2}  }

\newcommand{\balpha}{\boldsymbol{\alpha}}

\newcommand{\sgn}{\mathsf{sgn}}

\newcommand{\step}[1]{\eta_{#1}}

\newcommand{\qrec}{q_{\text{}}}

\newcommand{\vtheta}{\vec{\theta}}
\newcommand{\deterror}{\xi}
\newcommand{\vdeterror}{\vec{\deterror}}

\newcommand{\vzeta}{\vec{\zeta}}

\newcommand{\BoundarySet}{\cX_{\text{bound}}}
\newcommand{\Eboundary}{\mathcal{E}_{\text{boundary}}}
\newcommand{\Ebar}{\bar{\Ecal}}

\newcommand{\Mbar}{\bar{M}}

\newcommand{\hfunc}{h}

\newcommand{\Omegaconv}{\Omega_{\text{conv}}}

\newcommand{\vZ}{\vec{Z}}

\newcommand{\cX}{\mathcal{X}}
\newcommand{\cXbar}{\bar{\cX}}

\newcommand{\vPinit}{\vP_0}
\newcommand{\Treset}{\mathcal{T}}

\newcommand{\avproof}{\textcolor{black}}

\begin{document}
\maketitle
\begin{abstract}
We formulate the problem of fake news detection using distributed fact-checkers (agents) with unknown reliability. The stream of news/statements is modeled as an independent and identically distributed binary source (to represent true and false statements). 
Upon observing a news, agent $i$ labels the news as true or false which reflects the true validity of the statement with some probability $1-\pi_i$. In other words, agent $i$ misclassified each statement with error probability $\pi_i\in (0,1)$, where the parameter $\pi_i$ models the (un)trustworthiness of agent $i$. We present an algorithm to learn the unreliability parameters, resulting in a distributed fact-checking algorithm. Furthermore, we extensively analyze the discrete-time limit of our algorithm.
\end{abstract}
\section{Introduction}

As online social networks become increasingly effective in disseminating information, the task of distinguishing between true and false information becomes increasingly challenging. This growing efficiency of information dissemination has led to several studies on how misinformation spreads through networks \cite{acemoglu2010spread,acemoglu2021misinformation,budak2011limiting,papanastasiou2020fake,nguyen2012containment}. Conversely, there is growing interest in the development of automated fact-checkers that can perform tasks such as document retrieval, evidence extraction, and claim validation in an automated manner \cite{hanselowski2017framework,hanselowski2019richly,thorne2018fever}.
    
When there are multiple imperfect \textit{fact checkers}, determining the validity of a source based on their responses becomes a challenge. In such cases, it is important to know the reliability statistics of the {fact checkers} in question. As a result, a natural question arises: in the presence of multiple imperfect fact checkers, how can we formulate and learn their reliability over time?  We provide a model for distributed fact-checking using unreliable or imperfect agents. A key step in our model is to quantify the unreliability of each agent as the cross-over probability of certain BSC channels. {Given an estimate of the unreliability parameters, a weighted thresholding estimator can be used to identify the validity of the statement\cite{shapley1984optimizing,nitzan1981characterization,verma2023binary}, where the weights are the log-odds based on the agents' unreliability estimates.} We propose and study a learning rule to estimate the reliability parameters of the agents.  Our algorithm provides the advantage of requiring minimal memory and having a simplified update rule. 
     
In our problem, we are working with a mixture of product distributions. Determining the parameters of an identifiable mixture has been widely researched \cite{freund1999estimating,gordon2021source,gordon2022hadamard,chen2019beyond, feldman2008learning}. The parameter estimation problem typically involves finding a hypothetical model that produces samples with a distribution that closely resembles the true model.
The main contributions of the paper involve:
\begin{enumerate}[wide, labelwidth=!, labelindent=0pt, label=\roman*.]
    \item \textit{Formulation of Distributed Fact-Checking}: We study a model for distributed \textit{fact checking}, introduced in \cite{verma2023binary}, which constitutes agents modeled as Binary Symmetric Channels with unknown reliability in Problem~\ref{sec:problem_form}. 
    \item \textit{Online Estimator}: In Section~\ref{sec:estimator} we propose an online estimator for the unreliability parameters of the agents which makes use of the likelihood ratio between the source being fake or true given the agents' conclusion about the validity of the statement computed using the error estimate at a given time.
    \item \textit{Convergence Analysis}:  We study the convergence properties of the proposed online estimator involving truncated iterations.  To establish convergence, we utilize results akin to the Stochastic Approximation theorem presented in \cite{andrieu2005stability}. However, since the hypotheses of the theorem, specifically the assumptions on the Lyapunov function, in \cite{andrieu2005stability} are not met in our case, we extend the result and provide a proof tailored to our specific problem.
\end{enumerate}   
    \avnew{
    \textbf{Organization/Contribution of the Paper}: 
    After a brief review of the literature related to this topic,
    in Section~\ref{sec:problem_form}, we formally define the distributed fact-checking framework and describe the problem of interest to this work, namely the estimation of unreliability parameters, which is framed in Problem~\ref{prob:as_conv}. Additionally, we revisit the results from \cite{verma2024acc} that facilitate veracity estimation of news based on known unreliability parameters. Section~\ref{sec:estimator} introduces the specific estimator under consideration, including a necessary variant that incorporates projection/truncation, and elaborates on the intuition behind the estimator. 
    Section~\ref{sec:Main_Result} presents our main finding, Theorem~\ref{thm:discrete_time_convergence}, solving Problem~\ref{prob:as_conv} by showing convergence to a relevant ``equilibrium set''. Section~\ref{sec:Proof} presents the proof of Theorem~\ref{thm:discrete_time_convergence} by exploring the connection of the estimator to Stochastic Approximation, identifying a Lyapunov function, and studying the properties of the Lyapunov function that facilitate the convergence result. 
    }

    \textbf{Related Work}: 
        Given the unreliability parameters of the agents, an optimal approach to reconstructing unknown labels involves employing weighted majority voting. In this approach, the weights assigned to the output provided by each agent are equal to the log-odds based on the knowledge of the workers' unreliability~\cite{shapley1984optimizing,nitzan1981characterization}. In~\cite{verma2023binary}, we provide the characterization of weights that would result in the optimal estimator for labeling the validity of statements.
        
        On the other hand, the estimation of unreliability parameters is closely intertwined with the literature on crowdsourcing labeling where data labeling is crowdsourced to multiple unreliable workers. This process is susceptible to errors arising from various factors, such as task complexity, low incentive to accurately label the tasks, and the repetitive nature of the tasks. Estimating the unreliability of workers is challenging, as the true labels of the data are unknown.
        
        The challenge of distributed fact-checking shares similarities with crowdsourcing labeling tasks, which researchers have extensively studied within the Dawid-Skene model, introduced through empirical studies in 1979~\cite{dawid1979maximum}. The model emerged from medical applications where multiple clinicians label a patient's state. In this context, Dawid and Skene proposed an Expectation-Maximization (EM) algorithm. Over the years, various extensions and variants of this algorithm have emerged~\cite{raykar2010learning,albert2004cautionary,smyth1994inferring,hui1980estimating}, with a notable line of work employing spectral analysis of matrices representing correlations between agents and labeling tasks~\cite{zhang2016spectral}.
        
        Recent years have witnessed a growing body of research focused on performance guarantees for EM and its variants. Notably, Chao and Dengyong~\cite{gao2016exact}, as well as Zhang et al.~\cite{zhang2016spectral}, have provided performance guarantees for different versions of EM employing diverse initialization techniques.
        
        The convergence analysis of these variants of the Dawid-Skene estimator, rooted in the EM algorithm, has been explored for the offline scenario. In this context, where the sequence of statements to be verified is available as a batch, studies by Gao et al.~\cite{gao2016exact} and Zhang et al.~\cite{zhang2016spectral} have delved into the convergence aspects.
        
        The analyses of the EM-based algorithms hinge on a sufficiently accurate initialization derived from the output of a substantial batch of statements being validated. Importantly, all these works assume access to the storage of all labels of all agents, given their focus on an offline setting. The only notable work presenting an algorithm in a streaming setting, without the necessity to store the entire dataset, is found in the work of Bonald and Combes~\cite{bonald2017minimax}. Their proposed Triangular Estimation (TE) algorithm focuses on estimating the unreliability parameters of agents based on correlations between triplets of agents. This algorithm directly utilizes three agents, rather than the entire set, for estimating the unreliability parameter of a specific agent. The knowledge of all agents' output becomes indirectly relevant in determining which three agents to select for computing the unreliability parameter of a given agent. Our work is the first attempt at providing an online estimator which has similarities to the EM variants. In establishing convergence results for our online estimator, we draw connections to stochastic approximation concepts within the literature of control theory~\cite{andrieu2005stability}.
        
        \subsection{Notation}
        Let $\N$ denote the set of all natural numbers, $\N_0$ denote ${\N\cup \{0\}}$, and for any $n\in \N$, define ${[n]:= \{1,2,\dots,n\}}$. For any $i\in [n]$, let $[n]_{-i} := [n]\setminus\{i\}$. 
        We denote the set of real numbers by $\R$ and
        the set of all real-valued $n$-dimensional vectors by $\R^n$.
        Moreover, we use $\cX$ to denote the $n$-dimensional open unit cube, i.e., $\cX := (0,1)^n$.
        
        We use bold letters, such as $\vec{x}, \vec{s}$, to denote vectors, and regular letters, such as $x, s$, to denote scalars. For a scalar $a\in [0,1]$, we use $\bar a$ to denote $1-a$. We use $\one, \vec{0}$ to denote all-one and all-zero vectors, respectively, whose dimension will be clear from context. For a vector $\vec{x} \in \R^n$, $x_i$ denotes the $i$th element of $\vec{x}$. For any set $A$, we define ${d(\vx,A) := \inf_{\vec{y}\in A}  \norm{\vx-\vec{y}}}$ as the distance between $\vx$ and $A$. For any function $f:\mathcal{A} \to \mathcal{B}$ we define the function over the set $A \subseteq \mathcal{A}$ as $f(A) := \{f(x) \mid x \in A\}$. We define the indicator function as $\indicator{x\in A}= 1$ if $x \in A$ and $0$ otherwise.
        
        Throughout this work, all random variables are defined with respect to an underlying probability space $(\Omega, \Fcal, \Pr)$. When the probability measure is defined through a parameter, say $\vx$, we denote the probability measure by specifying $\vx$ as $\Pr(\cdot;\vx)$. Moreover, when the parameter $\vx$ is not specified, the probability measure is defined through the true parameter (described in the problem formulation) $\vpi$.
        
        For a sequence of entities such as $\{\vP(t)\}$, we denote the entry at time $t$ by $\vP(t)$. However, for step-size sequences specifically, we denote the step-size at time $t$ by using subscripts such as $\eta_t, \xi_t$. Throughout the paper, we use logarithm with respect to base $e$ and denote it by $\log x$ for any $x \in (0,\infty)$.
For a sequence of entities such as $\{\vP(t)\}$, we denote the entry at time $t$ by $\vP(t)$. However, for step-size sequences, we indicate the time $t$ by using subscripts, such as $\eta_t, \xi_t$. Throughout the paper, we use logarithm with respect to base $e$ and denote it by $\log x$ for any $x \in (0,\infty)$.
For a sequence of entities such as $\{\vP(t)\}$, we denote the entry at time $t$ by $\vP(t)$. However, for step-size sequences specifically, we denote the step-size at time $t$ by using subscripts such as $\eta_t, \xi_t$. Throughout the paper, we use logarithm with respect to base $e$ and denote it by $\log x$ for any $x \in (0,\infty)$.
  
    \section{Problem Formulation}\label{sec:problem_form}
    Consider a source that streams a sequence of statements where each statement can be true or false. We use a \emph{hidden} variable $S(t)\in \{+1,-1\}$ to denote the label (true/false) of the statement at discrete-time instance $t\in \N$. 
    We assume that the stream symbols are independently and identically distributed  according to the Rademacher distribution, i.e., 
    ${\Pr(S(t) = +1)  = \Pr(S(t) = -1)  = \frac{1}{2}}$, 
    for every $t \in \N$. A \emph{fact-checker} is interested in evaluating the validity of the statements using imperfect (inexpert) agents. 
    
    \textbf{Model for the fact checker:} 
    We model a \textit{fact-checker} as an overseer of multiple agents, where each agent is responsible for testing the validity of the statement provided to it. For $n\in \N$, let $[n]$ be the set of agents verifying the validity of the statements. At each time $t \in \N$, the agents provide noisy/imperfect labels or judgment for $S(t)$ to the fact checker by returning their assessment of the statement.
    In other words, if the agent considers the statement correct, it marks the statement as True. Otherwise, it marks it as False. However, due to their limited expertise, the agents' assessments may be different from the actual label of the statements.  Mathematically, we model \soh{the assessment at} agent $i\in[n]$ as a memoryless Binary Symmetric Channel (BSC) with the \textit{error probability} or \textit{crossover probability} $\pi_i \in (0,1)$, which takes the input ${S(t)}$ and outputs \soh{$R_i(t)$}, where for every $s\in\{\pm 1\}$, the distribution of the output is 
    \soh{
    $${\Pr(R_i(t) \hspace{-1pt}= -s|S(t)\hspace{-1pt}\!=\hspace{-1pt}s) \hspace{-1pt}=\hspace{-1pt} 1\!\!-\!\!\Pr(R_i(t) = s|S(t)\hspace{-1pt}=\hspace{-1pt}s ) = \pi_i}.$$
    }
    Therefore, agent $i\in [n]$ produces 
    an output $R_i(t)$, which is independent of the past.  Here, $\pi_i$ represents the unreliability of agent $i$ since the agent misclassifies the statement with probability $\pi_i$.
    Note that $\pi_i \in (0,1/2)$ embodies the fact that the agents are reliable on `average'. We represent the collection of crossover probabilities by $\boldsymbol{\pi}$ and the sequence of all agents' outputs at time $t$ by $\vec{R}(t)$.

    \textbf{Properties of Output distribution:} 
    Let us discuss some properties of output distribution. 
    \begin{enumerate}[wide, labelwidth=!, labelindent=5pt,label=\roman*.] 
        \item Since the statement stream $\{S(t)\}$ is independent, and each agent is viewed as a memoryless channel, the random vector process $\{\vec{R}(t)\}$ is an independent process. 
        \item At any time $t\in \N$, given $S(t)$, the outputs $\{R_i(t)\}_{i=1}^{n}$ are independent of each other. Moreover,  
        for any ${t\in \N}$ and for every ${i\in [n]}$, $R_i(t)$ has the Rademacher distribution.
        \item The joint distribution of the output $\vec{R}(t)$ given the true crossover parameters $\boldsymbol{\pi}$ is given as 
            \begin{align*}
          \hspace{-0.3cm}      \Pr(\vec{R}(t) \hspace{-1.6pt}=\hspace{-1pt} \vec{r};\vec{\pi}) &=\frac{1}{2}\hspace{-1pt}\left(\prod_{i=1}^n \pi_i^{\frac{1+r_i}{2}} \bpi_i^{\frac{1-r_i}{2}}  \hspace{-1pt}+\hspace{-1pt} \prod_{i=1}^n \pi_i^{\frac{1-r_i}{2}} \bpi_i^{\frac{1+r_i}{2}}  \hspace{-1pt}\right)\hspace{-1pt},      
            \end{align*}
        where $\vec{r}\in \{+1,-1\}^n$, and \avnew{ recall that} $\bar{x}=1-x$. 
    \end{enumerate}
    For brevity, given unreliability (crossover) parameters of the agents are $\vec{x}\in \cX$, we define $g_{\vec{x}}:\{+1,-1\}^n\to (0,1)$ to be the  distribution of the output vector $\vR \in \{-1,+1\}^n$, i.e., 
    \begin{align*}
        g_{\vec{x}}(\vec{r})\hspace{-1pt} &= \hspace{-1pt}\Pr(\vR \hspace{-1pt}=\hspace{-1pt} \vec{r};\vec{x}) =\frac{1}{2}\!\hspace{-1pt}\left(\prod_{i=1}^n x_i^{\frac{1+r_i}{2}} \bar{x}_i^{\frac{1-r_i}{2}}  \hspace{-2pt}+\hspace{-1pt} \prod_{i=1}^n x_i^{\frac{1-r_i}{2}} \bar{x}_i^{\frac{1+r_i}{2}}  \hspace{-2pt}\right)\hspace{-1pt}.
    \end{align*}
    In this notation, $g_{\boldsymbol{\pi}}(\vec{R}(t))$ refers to the true distribution of the output vector $\vec{R}(t)$ at any time $t\in\N$. 
    
   \av{
    The distributed fact-checker faces two key challenges. The first challenge is to \soh{estimate the true label of each statement. In other words, it should develop an}  estimator that minimizes the error probability in determining the validity of statements. To address this, \cite{verma2023binary} examined a class of estimators that make decisions based on whether linear combinations of individual agents' outputs exceed a certain threshold. \soh{Such an estimator can be expressed as}
    \begin{align}\label{def:linear_threshold}
         S_{\balpha,\tau}(\vec{R}) := \sgn \paran{\sum_{i=1}^n \alpha_i R_i-\tau} ,   
     \end{align}
     where $\vec{\alpha}\in \R^n$ and $\tau\in \R$ are given parameters. 
    In \cite{verma2023binary}, the set of optimal parameters for minimizing error probability $\Pr(S_{\balpha,\tau}(\vec{R}) \neq S)$ is characterized. It turns out that when dealing with a uniform\soh{ly distributed}  source (e.g., this work),  $S_{\balpha,0}(\vec{R})$ would be an optimal estimator where ${\balpha := \paran{\lpi{1},\dots, \lpi{n}} = \paran{\log \frac{1-\pi_1}{\pi_1},\dots,\log \frac{1-\pi_n}{\pi_n}}}$. This leads to the second challenge \soh{of the fact-checker, which is also}, the focus of this paper: estimating the agents' unreliability parameters $\vec{\pi}$.}
     
     In this paper, our goal is to obtain reliable 
     estimate for the unreliability parameters $\vec{\pi}$.  
     If the true statement validity sequence $\{S(t)\}$ were known, $\vec{\pi}$ could be easily estimated as the fraction of time at which $R_i(t)\neq S(t)$. The challenge here is to estimate the channel parameters without the knowledge of channel input. 
     \av{The ideal problem would be to identify an estimator that converges \emph{almost surely} (a.s.) to the true estimates of the unreliability parameter. However, some \soh{collections of parameters may result in indistinguishable distribution for the agents' output. It can be argued that those parameters cannot be distinguished from each other.}
     Thus, we focus on the following problem for multi-agent fact-checker. 
     \begin{prob}\label{prob:as_conv}
     Consider a fact-checker with access to the sequence $\{\vec{R}(t)\}$ of the assessments of $n$ agents, with unknown unreliability parameters $\pi_i$, for $i \in [n]$. 
     Determine an online estimator $\{\vP(t)\}$ of $\vec{\pi}$ based on the output of the $n$  agents such that 
     $\lim_{t\to \infty} d(\vP(t), \Scal) = 0$ a.s.,
     where $\Scal $ is the set of parameters that result in indistinguishable distribution for the agents' output, i.e., 
      \begin{align}\label{eq:S_def}
         \Scal := \{\vx \in \cX \mid g_{\vx}(\vr) = g_{\vpi}(\vr) \text{ for all } \vr \in \{\pm 1\}^n \}.
     \end{align}
     \end{prob}
     In fact, we can characterize the set $\Scal$ as follows.
     \begin{lem}\label{lem:Scharacterize}
         For $n\geq 3$ and 
          $\pi_i \in (0,1)\setminus\{{1}/{2}\}$ for $i\in [n]$, the set $\Scal$ defined in~\eqref{eq:S_def} is given by ${\Scal = \{\vpi, \one - \vpi\}}$.
     \end{lem}
     The proof of Lemma~\ref{lem:Scharacterize} is provided in  Appendix~\ref{appendix:Schar}.
     }
    \section{Estimator}\label{sec:estimator}
    First, let us introduce an online estimator for the unreliability parameters of the agents comprising the fact checker for any number of agents $n\geq 2$. 
    We have provided convergence guarantees for this algorithm for $n=2$ agents in~\cite{verma2024acc}.
   Consider the stream of output observed by the fact-checker $\{\vec{R}(t)\}$. \soh{At any time $t\in \N$, the fact-checker has an estimate $\vP(t)$ of the unreliability parameters $\boldsymbol{\pi}$, which depends on $(\vec{R}(1),\dots, \vec{R}(t))$. At time instance $t+1$ and upon observing the output $\vec{R}(t+1)$, the estimate can be updated to $\vP(t+1)$.} 
    Recall that if $\boldsymbol{\pi}$ was known, the fact checker could evaluate the likelihood ratio $ \frac{\Pr(\vec{R}(t+1)|S(t+1)=-1)}{\Pr(\vec{R}(t+1) |S(t+1)=+1)}$ to decode $S(t)$. Now, without $\boldsymbol{\pi}$, we can use its estimate $\vP(t)$ to compute an approximate likelihood ratio
    $L(t)$ of $S(t+1)=-1$ to $S(t+1)=+1$ based on $\vec{R}(t+1)$.  
    For this, let us define \soh{$L:\{+1,-1\}^n \times \cX\to \R$} by 
    \begin{align}\label{eqn:likelihood}
        L(\vR,\vx)=\prod_{i=1}^n \left(\frac{x_i}{1-x_i}\right)^{R_i}.
    \end{align}
    This represents the likelihood function of receiving $\vR$ given the unreliability parameters $\vx$.
    
    Now, for the received vector $\vec{R}(t+1)$ and an estimate $\vP(t)$ of their unreliability parameters, for brevity, let 
    \begin{align}\label{eqn:def_L}
        L(t) &:=L(\vR(t+1),\vP(t)).
    \end{align}
   Using $L(t)$, 
   we can estimate $S(t+1)$ by setting 
   \begin{align}
       \hat{S}(t+1) :=2\indicator{L(t)<1}-1 = \left\{
        \begin{array}{ll}
            -1 & \textrm{if $L(t)\geq 1$} \\
            +1 & \textrm{if $L(t)< 1$} 
        \end{array}
       \right. .
       \hspace{-5pt}
   \end{align}
    
    We are ready to discuss the update rule for the unreliability parameters' estimates, given the source symbol estimate ${\hat{S}(t+1)}$ and the output vector $\vec{R}(t+1)$. 
    Note that ${R_i(t+1)}$ agreeing with ${\hat{S}(t+1)}$ suggests that it is unlikely that the agent was introducing error at time $t$ and hence, we average $P_i(t)$ with a value less than half to obtain $P_i(t+1)$. 
    Similarly, if $R_i(t+1)$ disagrees with $\hat{S}(t+1)$, we average it with a value greater than half. 
    More precisely, the proposed algorithm/dynamics updates the unreliability parameters as
        \begin{align*}
            P_i(t+1) &= (1\hspace{-1pt}-\hspace{-1pt}\step{t}) P_i(t) \hspace{-1pt}+  \frac{1}{2}\step{t} \left( 1\hspace{-1pt}+\hspace{-1pt}\frac{L(t)- 1}{L(t)+1}R_i(t\hspace{-1pt}+\hspace{-1pt}1)\hspace{-1pt}\right)\!,
        \end{align*}
        for all $t\in \N$ and $i\in [n]$, with some initial condition ${\vP(0) \in (0,1)^n}$, 
        where $\{\step{t}\}$ is a pre-decided step-size sequence, and $L(t) $ is given in~\eqref{eqn:def_L}.        {Lext us rewrite the above iteration in compact form,\begin{align}\label{eq:update_rule}
            \vec{P}(t+1) &= \vec{P}(t) +\step{t}\ftilde(\vR(t+1),\vP(t)),
        \end{align}
        where ${\vftilde:\{+1,-1\}^n\times \cX\to \R^n}$ is the vector field with its $i$th coordinate \soh{is} given by \begin{align}\label{eq:ftildei_def}
            \ftilde_i(\vR,\vx) := \frac{1}{2}\paran{1+\frac{L(\vR,\vx)-1}{L(\vR,\vx)+1}R_i} - x_i.
        \end{align}
        
        Note that~\eqref{eq:update_rule} is a stochastic approximation-type iteration whose asymptotic behavior resembles the asymptotic behavior of the mean-field Ordinary Differential Equations (ODE)
        \begin{align}\label{eqn:MFODE}
            \dot \vx = \vf(\vx),
        \end{align}       
        where $\vec{f}\hspace{-1pt}:\hspace{-1pt}\cX \hspace{-1pt}\to\hspace{-1pt} \cX$ is defined by  
        \begin{align}\label{eq:fi_def}
            \vec{f}(\vec{x})
             &:= \E_{\vec{R}\sim g_{\vec{\pi}}}[\vftilde(\vR,\vx)].
        \end{align} 
        We will expand on this viewpoint in Section~\ref{sec:Proof}.

\subsection{Extension to Singly-extreme Vectors}        
        \beh{Note that the likelihood function~\eqref{eqn:likelihood} can be extended to the vectors $\vx\in[0,1]^n$, with \emph{only} one element being $0$ or $1$. Such a definition is not extendable to the case of vectors with more than \soh{one} extreme values.  \soh{For further elaboration,} let us first define \textit{singly-extreme vectors} as follows.  \begin{defn}\label{def:singly_extreme_points}
            For $i\in [n]$, let us define 
            \begin{align*}
            \BoundarySet^{(i)}\hspace{-1pt}:=\hspace{-1pt}\{\vx\in [0,1]^n\mid \hspace{-1pt}x_i \hspace{-1pt}\in\hspace{-1pt} \{0,1\},x_j\in(0,1)\forall j\in[n]_{-i}\}.
            \end{align*}
            We also define the set of all \emph{singly-extreme vectors } as ${\BoundarySet :=\bigcup_{i=1}^n  \BoundarySet^{(i)}}$.
             Furthermore, we define the extension of $\cX$ to include the set of singly-extreme vectors and denote it by
             ${\cXbar:= \cX \cup \BoundarySet}$. 
        \end{defn}}
        In Figure~\ref{fig:extreme_points_representation} (Left), we depict the sets $\cX, \BoundarySet^{1}$, and $\BoundarySet^{2}$ for the case of  $n=2$. 

        Assuming the convention ${\frac{1}{0}:=\lim_{p\to 0^+}\frac{1-p}{p}=\infty}$,  for a singly-extreme vector ${\vx\in \BoundarySet^{(i)}}$, we define 
        \begin{align*}
            L(\vR,\vx):=\begin{cases}
                0&\text{if $(-1)^{x_i}=R_i$}\\
                \infty&\text{if $(-1)^{x_i}\not=R_i$}
            \end{cases},
        \end{align*}
        leading to 
        \begin{align}\label{eqn:LratioExtension}
            \frac{L(\vR,\vx)-1}{L(\vR,\vx)+1}=\begin{cases}
                -1&\text{if $(-1)^{x_i}=R_i$},\\
                +1&\text{if $(-1)^{x_i}\not=R_i$}.
            \end{cases}
        \end{align}
        Note that if $x_i,x_j\in \{0,1\}$ for $i\not=j$,  then the likelihood ratio~\eqref{eqn:likelihood} cannot be defined for all vectors $\vR\in\{+1,-1\}^n$. In particular, if $(-1)^{x_i}R_i\not=(-1)^{x_j}R_j$, the product in~\eqref{eqn:likelihood} would contain a $0$ and $\infty$ term leading to an undefined expression ${0\times \infty}$. This is, in fact, a fundamentally unresolvable phenomenon as this is related to the case where the fact-checker is receiving two contradictory verdicts for the same statement from two fully reliable/unreliable agents. 

        With this discussion in mind, we can extend the definition of $\ftilde$ in \eqref{eq:ftildei_def} and $\vf$ in \eqref{eqn:MFODE} to the set $\cXbar= \cX \cup \BoundarySet$, by considering the ratio~\eqref{eqn:LratioExtension} for $\vx\in \BoundarySet^{(i)}$ for $i\in[n]$.
    \begin{figure}
    \centering
    \begin{tikzpicture}[scale=4] 
        \coordinate (A) at (0,0);
        \coordinate (B) at (1,0);
        \coordinate (C) at (1,1);
        \coordinate (D) at (0,1);
    
        \draw[ao, line width=1.5pt] (A) -- (B);
        \draw[ao, line width=1.5pt] (C) -- (D);
        
        \fill[pattern=north east lines] (A) rectangle (C); 
    
        \draw[red, line width=1.5pt] (A) -- (D);
        \draw[red, line width=1.5pt] (B) -- (C);
    
        \foreach \point in {A,B,C,D}
            \fill[blue] (\point) circle (0.05); 
    \end{tikzpicture}\hspace{1cm}
     \includegraphics[width=0.3\linewidth]{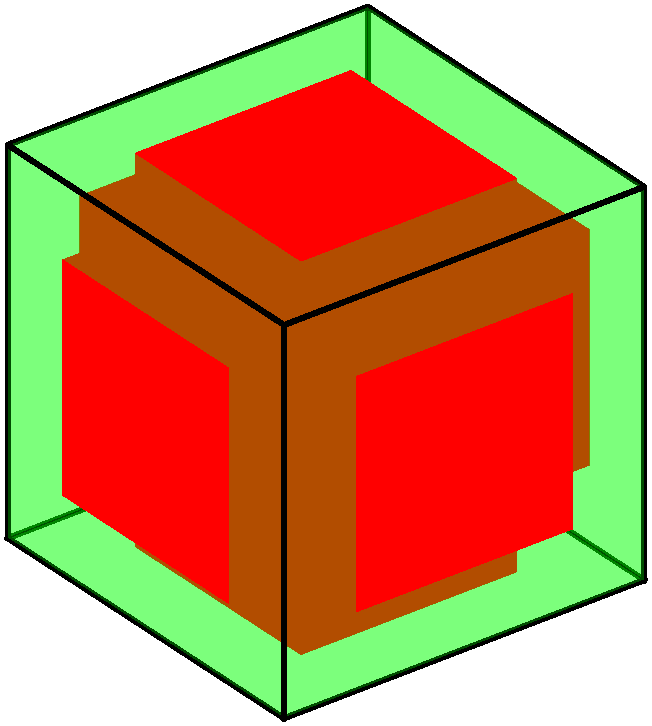}
    \caption{Left: Parameter set for $n=2$ agents: the red and green lines represent the sets $\BoundarySet^{(1)}$ and $\BoundarySet^{(2)}$ respectively. The shaded region represents $\cX$. The box excluding the blue points represents $\cXbar$. Right: Truncation set for $n=3$ agents: the red region represents the form of a set $\Kcal_t$. 
    The cube excluding the solid lines show $\cXbar$.}
    \label{fig:extreme_points_representation}
\end{figure}    
\subsection{The Step Size}
        For the step-sizes $\step{t} $, we assume they satisfy the following stochastic approximation step-size assumption. 

        \begin{assum}\label{asmp:stepsizes}
        The step-sizes $\{\step{t}\}$ are positive, non-increasing, and  satisfying $\sum_{t=0}^\infty \step{t} =\infty$ and ${\sum_{t=0}^\infty \step{t}^2 <\infty}$.
        \end{assum}

        One popular choice for the step-size sequence is the harmonic sequence $\step{t} = \frac{1}{t+1}$ for all ${t\in \N}$. 
        {To grasp the motivation behind the estimator using such a step-size sequence, examine the scenario when the fact checker \textit{knows} the source sequence symbols $\{S(t)\}$. 
        Since, at any time $t\in \N$, conditioned on $S(t)$, the output distributions of the agents \soh{are} independent of each other, the problem of estimating the unreliability parameters of the agents is equivalent to $n$ uncoupled problems of estimating the parameter of a Bernoulli distribution from its independent samples. The estimation of the parameters for this problem is well-studied. \soh{In particular, an effective class of estimators to solve this problem} is the \textit{add-constant} estimator~\cite{kamath2015learning}. For the current setting, for each $i\in [n]$ the add-$\beta$ estimator, where $\beta \geq 0$, for parameter $\pi_i$ at time $t\in \N$ is given by 
        \begin{align*}
            Q_i(t) = \frac{\beta+\sum_{k=1}^{t} \indicator{R_i(k)\neq S(k)}}{t+2\beta}. \nonumber
        \end{align*}
        The estimator makes use of the empirical frequency of agent~$i$ misclassifying the source symbol  and can be expressed recursively as 
        \begin{align*}
            Q_i(t+1) = (1-\stepq_{t}) Q_i(t) + \stepq_{t} \indicator{R_i(t+1)\neq S(t+1)}.
        \end{align*}
        Here, $\stepq_{t}:= \frac{1}{t+1+2\beta}$ and $Q_i(0)=1/2$.
        The convergence properties of estimator $Q(t)$ for different values of $\beta$ and various loss functions are studied in~\cite{kamath2015learning}. Different values of $\beta$ lead to well-known estimators, including the empirical estimator \mbox{($\beta=0$)},  the Krichevsky–Trofimov (KT) estimator ($\beta=\frac{1}{2}$), and Laplace estimator ($\beta=1$).}

        To see the connection to our setting, where the source symbol is unknown, consider an extreme case where $L(t)\gg 1$ (which implies $\hat{S}(t+1)=-1$). For $R_i(t+1)=+1$, we get 
        \[\frac{1}{2}\bigg(\frac{L(t)- 1}{L(t)+1}R_i(t+1) +1\bigg) = \frac{L(t)}{L(t)+1}\approx 1,\] whereas for $R_i(t+1)=-1$  we have  ${\frac{1}{L(t)+1}\approx 0}$. Thus, 
        \[
        \frac{1}{2}\bigg(\frac{L(t)- 1}{L(t)+1}R_i(t+1) +1\bigg) \approx \indicator{R_i(t+1)\neq \hat{S}(t+1)}.
        \]
        A similar situation holds when $L(t)\approx 0$. 
        Therefore, the update rule~\eqref{eq:update_rule} with $\step{t}=\frac{1}{t+1}$ can be viewed as an imperfect and adaptive version of the add-$\beta$ estimator with $\beta=0$.


\subsection{The Truncation Sets}

    With the presence of multiple agents making decisions, in order to obtain convergence guarantees, it is important to stay away from the boundary where two or more agents have estimated unreliability close to $0$ or $1$. 
    In particular, in line with the concept of expanding truncations for stochastic approximation in \cite[Chapter 2]{chen2005stochastic}, we \soh{introduce} a collection of growing truncation sets. These compact sets gradually converge to a superset. The estimate (but not the process) is reset to an arbitrary initial value each time it crosses the \emph{current} truncation set.

    We define $\{\Kcal_t\}$ to be a sequence of increasing compact truncation sets for $\cXbar$, i.e., (i) ${\bigcup_{t=0}^{\infty} \Kcal_t = \cXbar}$, (ii)~$\Kcal_t$ is compact for all $t \in \N_0$, and (iii) $\Kcal_t \subseteq \Kcal_{t+1}$ for every $t \in \N_0$. 
    There are many choices for the sequence of increasing truncation sets $\Kcal_t$ that satisfy \soh{these three conditions}. 
One such class of sets has the form 
${\Kcal_t \!=\! \bigcup_{i=1}^n \{\vx \hspace{-1pt}\in\hspace{-1pt} [0,1]^n \!\! \mid \!\! x_i\in
[0,1],
|x_j\hspace{-1pt}-\hspace{-1pt}\frac{1}{2}|\leq r_t, \forall j \hspace{-1pt}\in\hspace{-1pt} [n]_{-i} \},}$
where $\{r_t\}$ is a strictly increasing sequence satisfying ${\lim_{t\to\infty} r_t = \frac{1}{2}}$. 

\subsection{Estimator with Resettings}
Consider the set of increasing truncation covers $\{\Kcal_t\}$ and an initial estimate $\vPpr(0)\in \Kcal_0$. 
Let $\{\proj{t}\}$ be a sequence of non-negative integers that keeps track of the current \soh{truncation} set of the algorithm with $\proj{0}=0$. Then recursively, for any $t \in \N_0$, consider the update of the estimate \beh{similar} to~\eqref{eq:update_rule} at point $\vPpr(t)$, \soh{that is,}  
${\vec{y} = \vPpr(t) + \step{t} \vftilde(\vR(t+1),\vPpr(t))}$. We reset the dynamics if $\vec{y}$ is outside the current active set $\Kcal_{\proj{t}}$, otherwise, let $\vPpr(t+1)=y$. In other words, 
\begin{align}\label{eq:projected_update_rule}
       \vPpr(t+1) = 
       \begin{cases}
           \vec{y},   &\text{ if } \vec{y} \in \Kcal_{\proj{t}} \\
        \vPinit,        &\text{ if } \vec{y} \not \in \Kcal_{\proj{t}}
       \end{cases},
    \end{align}
where
$\vPinit\in \Kcal_0$ 
is an arbitrarily chosen point from the set $\Kcal_0$, and the counter for the truncation set is defined as 
\begin{align}\label{eq:counter_def}
    \proj{t+1} = \proj{t} + \indicator{\vec{y} \not \in \Kcal_{\proj{t}}}. 
    \vspace{-5pt}
\end{align}
\vspace{-5pt}
\section{Main Results}\label{sec:Main_Result}
In the next theorem, we 
offer a solution to Problem~\ref{prob:as_conv} by showing that the online estimator defined in~\eqref{eq:projected_update_rule}, almost surely converge to \soh{the} intended set. To articulate the theorem, define 
$${\hfunc(a,b):= ab+(1-a)(1-b)}$$ for ${a, b \in [0,1]}$.

\begin{thm}[Convergence of Online Estimator]\label{thm:discrete_time_convergence}
For a fact-checker comprised of agents with unreliability parameters ${\vec{\pi} \in \cX}$, under Assumption~\ref{asmp:stepsizes} on the step-sizes, with probability one the online estimator $\{\vPpr(t)\}$, defined in~\eqref{eq:update_rule}, converges to the set \beh{$ \Ebar =\Ecal \cup \Eboundary$
where 
$\Ecal$ is the set of equilibrium points of the mean-field ODE~\eqref{eq:fi_def}, i.e., $\Ecal = \{\vec{x}\in \cX\mid \vf(\vx)=\vec{0}\}$ 
and 
$\Eboundary$ is the collection of $2n$ points, \soh{given by}
 \begin{align}\label{eq:Eboundary_def}
         \!\Eboundary := \bigcup_{i=1}^n\{ &\vx \in \BoundarySet\mid \,  x_i\in \{0,1 \}, 
        &\hspace{-2.5pt}x_j\hspace{-1pt}=\hspace{-1pt}\bar{x}_i h(\pi_i,\bar{\pi}_j) \hspace{-1pt}+\hspace{-1pt} {x_i} h(\pi_i,{\pi}_j)
        , \forall j \hspace{-1pt}\in\hspace{-1pt} [n]_{-i}\}. 
    \end{align}}  
\end{thm}

To understand the set $\Ebar$ to which the estimates converge, note that it is comprised of $\Ecal$, \av{the set of zeros of the mean-field ODE $\vf(\vx)=\E_{\vR\sim g_{\vec{\pi}}}[\vftilde(\vR,\vx)]$}
and $2n$ points at the boundary of the set $\cXbar$.  
We know that points $\vpi$ and $\one - \vpi$ are present in the set $\Ecal$ and they should naturally be there as $\vpi$ and $\one - \vpi$ are indistinguishable from the distribution of the output vector $\vR$. \av{Moreover, we conjecture that the set $\Ecal$ is comprised of these two points and the trivial point $\frac{1}{2}\one$. 
\begin{conj}[Characterization of $\Ecal$]\label{conj:ecal}
    For $n\geq 3$ and $\pi_i \in (0,1)\setminus\{1/2\}$ for $i\in[n]$, we conjecture that $\Ecal = \{\vpi, \one - \vpi, \frac{1}{2}\one\}$.
\end{conj}
}
\avnew{For $n=3$, the  conjecture above has been established in \cite{verma2024cdc}.}

On the other hand, to understand the points on the boundary, note that $h(\pi_i,\bar{\pi}_j)$ represents the probability of agents $i$ and~$j$ declaring the opposite verdicts regarding the validity of the statement, \soh{i.e., $ h(\pi_i,\bar{\pi}_j) = \Pr(R_i \neq R_j) $, and similarly, $ h(\pi_i,{\pi}_j) = \Pr(R_i = R_j)$. The estimate $x_i = 0$ implies that  the fact-checker has decided that agent~$i$ is completely reliable and hence $S=R_i$; then, for agent~$j$, instead of estimating ${\pi_j= \Pr(R_j \neq S)}$, it is estimating ${x_j = \Pr(R_j \neq R_i)= h(\pi_i,\bar{\pi}_j)}$. Similarly, if ${x_i=1}$, the fact-checker has decided that agent $i$ is completely unreliable, i.e.,  ${R_i=-S}$; hence, for agent $j$, the fact-checker is estimating ${x_j = \Pr(R_j \neq S)= \Pr(R_j = R_i)  = h(\pi_i,\pi_j)}$}.
\section{Proof of Main Theorem}\label{sec:Proof}
In this section, we provide the proof of our main result, Theorem~\ref{thm:discrete_time_convergence}. The road map to the proof is structured as follows:

\avnew{
 \begin{enumerate}[wide, labelwidth=!, labelindent=5pt,label=(\arabic*)] 
        \item Lemma~\ref{lem:StochApprox} expresses the estimator in~\eqref{eq:update_rule} as a stochastic approximation of a mean-field ODE and characterizes the corresponding vector field $\vf$.  Lemma~\ref{lem:rewritten_f}  presents a simplified mean-field ODE that emphasizes its connection to the estimator of the label via $\inp{\vr}{\bl_{\vx}}$ (compare the presence of the expression in \eqref{def:linear_threshold}), for a specified unreliability parameter $\vx$). 
        \item  An essential element in establishing the proof of the main theorem involves identifying a Lyapunov candidate  (Theorem~\ref{thm:KL_lyapunov}) and studying the properties of its sublevel sets which will be established in Lemma~\ref{lem:V_asmp_satisfied}. In Lemma~\ref{lem:Vvalue_extremepoint}, we evaluate the mean-field function $\vf$ and the Lyapunov function $V$ at the \emph{singly-extreme points} $\cXbar$ (introduced in Definition~\ref{def:singly_extreme_points}). This assessment shows that for a sufficiently large level, the closure of the corresponding sublevel set of the Lyapunov function includes a subset of $\cXbar$.  Lemma~\ref{lem:V_asmp_satisfied} specifies the characteristics of the Lyapunov function $V$, which serves as a central element in the proof.
        \item 
        Lemma~\ref{lem:Kq_recurrence} utilizes the properties of the Lyapunov function's sublevel set, those of the truncation set, and the boundedness of the accumulated error shown in Lemma~\ref{lem:error_convergence}, to ensure that after a long enough time the updates remain within a compact set.
        For each sample path with bounded accumulated error (sample path corresponding to $\omega \in \Omegaconv$ as defined in \eqref{eq:omegaConv}) as the updates stay within a compact subset, we demonstrate in Theorem~\ref{thm:stability_sa_main} that the updates converge to a specified set. Ultimately, by integrating Lemmas~\ref{lem:error_convergence}, 
        ~\ref{lem:Kq_recurrence}, and Theorem~\ref{thm:stability_sa_main}, we provide the proof of our main result.
    \end{enumerate}
    }

\subsection{Stochastic Approximation}
\soh{A stochastic approximation of an Ordinary Differential Equation (ODE) is a recursive method to find/approximate the solutions of an ODE. For a mean-field ODE $\dot{\vec{z}}= \vec{F}(\vec{z})$ with  $\vec{F}:\mathcal{Z} \to \R^d$, consider the sequence $\{\vZ(t)\}_{t\geq 0}$, generated by 
\begin{align}\label{eq:intro_SA}
    \vZ(t\hspace{-1pt}+\hspace{-1pt}1) &\hspace{-1pt}=\hspace{-1pt} \vZ(t) + \step{t} (\vec{F}(\vZ(t)) + \vdeterror(t\hspace{-1pt}+\hspace{-1pt}1)),  
\end{align}
for $t\in \N_0$. Here, $\vZ(0)\in \mathcal{Z}$ is an initial point, and $\{\step{t}\}$ is a step-size sequence that is often assumed to satisfy Assumption~\ref{asmp:stepsizes}. Furthermore, \{$\vdeterror(t)\}$ is a martingale-difference sequence, i.e.,  
${\E [\vdeterror(t+1) \mid \Fcal_t] = 0}$
for all $t\in \N_0$, where ${\Fcal_t = \sigma(\vec{Z}(k),\xi(k):k\leq t)}$ is the $\sigma$-algebra generated by the past. Then, the update vectors $\vZ(t)$ can be viewed as an approximation of the path taken by the mean-field ODE. More precisely, if $\vec{F}$ is a Lipschitz function, $\vdeterror(t)$ is bounded in expectation, and $\{\vZ(t)\}$ remains bounded, then, the sequence $\{\vZ(t)\}$ converges to the solution of $\dot{\vec{z}}= \vec{F}(\vec{z})$~\cite{borkar2009stochastic}.} 

To prove the main result, we first show that the difference in the updates of the estimates in~\eqref{eq:update_rule} can be decomposed into a deterministic part and a zero-difference martingale. 
\begin{lem}\label{lem:StochApprox}
For $t\ \in \N_0$, the online estimator~\eqref{eq:update_rule} satisfies 
\begin{align*}
    \vP(t+1) = \vP(t) + \step{t} \left( \vec{f}(\vP(t)) + \vM(t+1) \right),
\end{align*}
where $\vec{f}$ is the mean-field ODE given in~\eqref{eq:fi_def} and 
\[
\vM(t+1):=\vftilde(\vR(t+1),\vP(t))-\vf(\vP(t)).
\]
Moreover, for the sequence $\{\vM(t)\}$, \soh{we have 
\begin{enumerate}[label=(\roman*)]
    \item $\E[\vM(t+1)\mid \Fcal_t] = 0$ and
    \item $\|\vM(t+1)\|_{\infty} \leq 2$ (and hence, bounded in-expectation),
\end{enumerate}
with probability one for all $t \geq 0$. Here, }the filtration $\{\Fcal_t\}$ is defined as ${\Fcal_t = \sigma(\vP(k), \vM(k): k\leq t)}$ for all $t \in \N_0$. 

In other words $\{\vM(t)\}$ is a bounded martingale difference sequence with respect to the filtration $\{\Fcal_t\}.$
\end{lem}
The proof of Lemma~\ref{lem:StochApprox} is provided in Appendix~\ref{appendix:stoch_approx}.

Although the discrete-time process $\{\vP(t)\}$ satisfies the conditions to be viewed as the Stochastic Approximation scheme for $\vec{f}(\vec{x})$, we cannot use the standard results in stochastic approximation \cite{borkar2009stochastic, kushner2006stochastic} \soh{to guarantee} the convergence of our learning rule,  since the function $\vf$ is not a Lipschitz function.
Instead, we use certain Lyapunov functions and their properties to show the convergence of the proposed online estimator. 

In the following lemma, we provide an alternative \soh{representation} for the probability of the output vector $g_{\vx}(\vr)$ and the update value $\vftilde(\vr,\vx)$. We refer to Appendix~\ref{appendix:lyapunov_ode_prop} for the proof. 
\begin{lem}\label{lem:rewritten_f}
For any $\vr \in \{-1,+1\}^n$ and $\vx \in \cX$, the expression of output probability for $\vR$ based on unreliability vector $\vx$, $\Pr(\vR=\vr;\vx)$, is given as 
\begin{align}\label{eq:output_prob_alternate}
    g_{\vx}(\vr)\!\!  &=\!\! \paran{\frac{1}{2} \exp {\frac{-\inp{\vr}{\bl_{\vx}}}{2}}\!  +\! \frac{1}{2} \exp {\frac{\inp{\vr}{\bl_{\vx}}}{2}} }
      \!\!\prod_{i=1}^n \sqrt{x_i(1\hspace{-1pt}-\hspace{-1pt}x_i)} \nonumber\\
      &= \cosh \paran{\frac{\inp{\vr}{\bl_{\vx}}}{2}} \prod_{i=1}^n \sqrt{x_i(1-x_i)}, 
\end{align}
\soh{where the $i$th coordinate of vector $\bl_{\vx}$ is given by  ${\bl_{x_i} = \log\frac{1-x_i}{x_i}}$.}
Moreover, the function $\vftilde(\vr,\vx)$ \soh{introduced in~\eqref{eq:ftildei_def}} can be expressed as 
\begin{align}
    \vftilde(\vr,\vx)   
        &= \frac{1}{2}
        \paran{
           \one -\tanh\paran{  \frac{1}{2} \inp{\vr}{\bl_{\vx}} } \vr
        } - \vx. \label{eq:rewritten_f1}
\end{align}
\end{lem}

\soh{Next, we introduce a correction term $\vcorrection{t}$ to express the resetting based update rule \eqref{eq:projected_update_rule} in a recursive stochastic approximation form. The correction term becomes active only at the resetting iterations. Later in Lemma~\ref{lem:Kq_recurrence}, we show that the correction term $\vcorrection{t}$ is non-zero only finitely often.}

\beh{The update rule $\{\vPpr(t)\}$ defined in \eqref{eq:projected_update_rule} can also be expressed in the stochastic approximation form as 
    \begin{align}\label{eq:Ppr_def_Mt}
        \vPpr(t\hspace{-1pt}+\hspace{-1pt}1)\hspace{-1pt} = \hspace{-1pt}\vPpr(t) \hspace{-1pt}+\hspace{-1pt} \step{t}[\vf(\vPpr(t)) \hspace{-1pt}+\hspace{-1pt} \vM({t\hspace{-1pt}+\hspace{-1pt}1})+ \vcorrection{t\hspace{-1pt}+\hspace{-1pt}1}],
    \end{align}
    where $\vf(\cdot), \vM(t+1)$ are as discussed in Lemma~\ref{lem:StochApprox}, and the correction term for the event of resetting is given as  \begin{align*}
        \vcorrection{t+1} := &\frac{1}{\step{t}}\indicator{\gamma(t+1)\not=\gamma(t)}\times 
         \left(\vPinit - (\vPpr(t) + \step{t} \vftilde(\vR(t+1),\vPpr(t))\right), \nonumber
    \end{align*}
    and $\vPinit\in \Kcal_0$ is the reset estimate.}
\subsection{Lyapunov Function}\label{sec:lyap_function}

In this section, we propose and study a Lyapunov function for the mean-field ODE \eqref{eqn:MFODE}, 
and we prove the desirable properties related to the function to establish the almost sure convergence of the update rule~\eqref{eq:projected_update_rule}. To define the Lyapunov function, we utilize the Kullback-Leibler (KL) divergence which is defined as follows.  
\begin{defn}[Kullback-Leibler Divergence \cite{polyanskiy2014lecture}]
Let $\mu, \nu$ be \soh{two} distributions on discrete set/alphabet $\mathcal{U}$. The KL-divergence between $\mu$ and $\nu$ is defined as 
${\Dkl(\mu \| \nu) = \sum_{u\in \mathcal{U}} \mu(u)\log \frac{\mu(u)}{\nu(u)}}$, where we use the conventions (i) $0\cdot\log \frac{0}{0} = 0$, (ii) if there exists $u\in \mathcal{U}$ such that $\nu(u)=0$ and $\mu(u)>0$ then $\Dkl(\mu \| \nu) = \infty.$
\end{defn}

For a Lyapunov candidate, we focus on the Kullback-Leibler divergence between the \soh{distribution of the output $\vR$ generated by the agents when the unreliability parameters vector is $\vec{\pi}$ or $\vx$}. More precisely, for $\vec{\pi} \in \cX$ we define the Lyapunov candidate as $V:\cX\to [0,\infty)$ with 
\begin{align}\label{eqn:lyap}
    {V(\vec{x}) := \Dkl(g_{\boldsymbol{\pi}}\| g_{\vec{x}})}.
\end{align} \beh{Note that here the alphabet is $\mathcal{U}=\{+1,-1\}^n$.}
In addition, the Lyapunov function is finite for all $\vx \in \cX=(0,1)^n$ as $g_{\vx}(\vr) > 0$ for every $\vr \in \{+1,-1\}^n$. \soh{Moreover, $V(\vx)$ can be rewritten as 
\begin{align*}
    V(\vx) =  C_{\vpi} - \sum\nolimits_{\vr \in \{+1,-1\}^n}  g_{\vpi}(\vr)\log {g_{\vx}(\vr)},
\end{align*}
where $C_{\vpi} := \sum_{\vr \in \{+1,-1\}^n}  g_{\vpi}(\vr)\log {g_{\vpi}(\vr)}$ is a finite constant for any $\vpi \in \cX$. Since the second term is a convex combination of $2^n$ functions that are logarithms of polynomials in ${\vx}$, we can conclude that $V(x)$ is
a continuously differentiable function.}

\textbf{Extension to $\cXbar$:} Note that the Lyapunov candidate function~\eqref{eqn:lyap} can be extended to the set of singly-extreme vectors~${\vx \in \BoundarySet}$.
For $r\in\{+1,-1\}$ and $x\in (0,1)$ we can rewrite 
$x^{\frac{1+r}{2}} (1-x)^{\frac{1-r}{2}} = \frac{1+r}{2}x+\frac{1-r}{2}(1-x)$.  
From the above representation for $x=0$, we adopt the convention $0^{\frac{1+r}{2}} := \lim_{x\to 0} \frac{1+r}{2}x+\frac{1-r}{2}(1-x) = \frac{1-r}{2}$ for $r\in \{-1,+1\}$. 
Based on the convention, we can extend the definition of the output distribution, $g_{\vx}$, for $\vx \in \BoundarySet$. For $i\in[n]$, for a singly-extreme vector $\vx \in \BoundarySet^{(i)}$ and vectors $\vr \in \{+1,-1\}^{n}$, the output distribution is defined as 
\begin{align}\label{eq:gx_extend}
    g_{\vx}(\vr) :=
    \begin{cases}
       \prod_{k\in [n]_{-i}} x_k^{\frac{1-r_k}{2}}(1-x_k)^{\frac{1+r_k}{2}} & \text{ if } r_i = (-1)^{x_i} \\
       \prod_{k\in [n]_{-i}}  x_k^{\frac{1+r_k}{2}}(1-x_k)^{\frac{1-r_k}{2}}  & \text{ if } r_i \neq (-1)^{x_i}
    \end{cases}.
\end{align}
For $\vx \in \BoundarySet$, we extend the definition of $V(\vx)$ using the definition in~\eqref{eq:gx_extend} for $g_{\vx}$.

Moreover, for $i\in[n]$, we define the gradient of $V(\vx)$, \soh{denoted by} $\nabla V(\vx)$,  for $\vx \in \BoundarySet^{(i)}$, through its $j$th element 
\begin{align*}
    \del{V(\vx)}{x_j} = 
    \begin{cases}
        \lim_{h \to 0^+} \frac{V(\vec{x}+(-1)^{x_i}h \vec{e}_i)-V(\vec{x})}{(-1)^{x_i}h} & \text{ if } j = i \\
        \lim_{h \to 0} \frac{V(\vec{x}+h \vec{e}_j)-V(\vec{x})}{h} & \text{ if } j \neq i 
    \end{cases},
\end{align*}
where $\vec{e}_i$ is the $i$th standard basis vector of $\R^n$. 

In the following theorem, we show that $V(\vx)$ satisfies the conditions to be a Lyapunov function for the ODE ${\dot \vx = \vf(\vx)}$. 

\begin{thm}\label{thm:KL_lyapunov}
For the mean-field dynamics~\eqref{eqn:MFODE} and the function $V(\vec{x})$ given in~\eqref{eqn:lyap}, we have $\inp{\nabla V(\vec{x})}{ \vec{f}(\vec{x})} \leq 0$ for all $\vec{x} \in \cX$. Furthermore the equality holds iff $\vec{x} \in \Ecal$, where ${\Ecal= \{\vec{x}\in \cX \mid \vec{f}(\vec{x}) = \vec{0}\}}$ is the set of equilibrium points of the ODE.  
\end{thm}
The proof of Theorem~\ref{thm:KL_lyapunov} is provided in Appendix~\ref{appendix:lyapunov_ode_prop}. 
%

To prove the convergence, we use the Stochastic Approximation result in~\cite{andrieu2005stability}. 
In~{\cite[Theorem 2.3]{andrieu2005stability}}, it is established that under certain assumptions \soh{on} the functions $\vf(\vx)$ and $V(\vx)$, if the updates stay in a compact subset $\Kcal$ of $\cX$, and the step-sizes and error term satisfy certain boundedness conditions, then the updates converge to the set $\Kcal \cap \Ecal$. 
In doing so a key property that is being used is the compactness of the sublevel set of Lyapunov function $V(\vx)$ associated with the ODE $\dot \vx = \vf(\vx)$. 
\av{However, note that the Lyapunov function $V(\vx) = \Dkl(g_{\vpi} \| g_{\vx})$ defines a non-compact sublevel set for levels greater than a certain value. 
In other words, for large enough $M$,  \soh{the set} 
${\{\vx \in \cX \mid V(\vx) \leq M\}}$ 
is not a compact set, as will be clarified as a consequence of Lemma~\ref{lem:Vvalue_extremepoint} in Section~\ref{sec:boundary_bhvr}.}
Consequently, if we focus only on the updates over $\cX,$ the compactness assumption essential for \soh{Theorem}~\ref{thm:stability_sa_main}, \avnew{that proves the convergence behavior of particular trajectories (linked to any $\omega \in \Omegaconv$) of the stochastic process $\{\vPpr(t)\}$}, would not hold. 
However, we can address this issue by extending the function to the set $\cXbar$ ensuring the assumptions are satisfied, as detailed in the following subsection.}
\subsection{Boundary Behavior with Extreme Unreliability}\label{sec:boundary_bhvr}

In this subsection, we discuss the behavior of the Lyapunov function $V(\vx)$ and the mean-field function $\vf(\vx)$ for the singly-extreme vectors (see Definition~\ref{def:singly_extreme_points}). To this end, we expand the definition of both functions for $\vx \in \BoundarySet$ as the limit taken along any trajectory inside $\cX$. With these values, we show that the Lyapunov function takes finite value over $\BoundarySet$. \av{Note that for $\vx \in \BoundarySet$, the definitions of the functions for $\vf(\vx)$ and $V(\vx)$ in \eqref{eq:fi_def} and \eqref{eqn:lyap} respectively remain valid and yield the expressions presented in Lemma~\ref{lem:Vvalue_extremepoint}.
} 

To introduce the next lemma, let's define $H_a(x)$ for ${a,x \in (0,1)}$ as
$H_a(x):= -a\log x - (1-a)\log(1-x).$ 
\begin{lem}\label{lem:Vvalue_extremepoint}
    Define $C_{\vpi} := \E_{\vR\sim g_{\vec{\pi}}}[\log g_{\vec{\pi}}(\vR)]$. Then for \beh{any $i\in[n]$ and $\vx \in \BoundarySet^{(i)}$}, we have
    \begin{enumerate}[label=(\roman*).]
        \item 
        ${V(\vx)  = C_{\vpi} + \log 2 + \sum_{k \in [n]_{-i}} H_{h(\pi_i,\pi_k)}(\hfunc(x_i,{x}_k)),}$
        \item $\vf(\vx)$ satisfies $f_i(\vx)= 0$ and 
        \begin{align*}
            f_j(\vx) = \hfunc(\pi_i,h(x_i,\pi_j))- x_j,  \quad \forall j \in [n]_{-i}.
        \end{align*}
        \item $f_j(\vx) = x_j(1-x_j)\del{V(\vx)}{x_j}$, for all $j \in [n]$. 
    \end{enumerate}
\end{lem}
The proof of Lemma~\ref{lem:Vvalue_extremepoint} is provided in Appendix~\ref{appendix:proof_extreme}.
\begin{cor}
    For $\vx \in \BoundarySet$, $\vf(\vx) = 0$ iff $\vx \in \Eboundary$ where $\Eboundary$ is given in \eqref{eq:Eboundary_def}.
\end{cor} 
Since the function $H_a(x)$ is minimized at $x = a$, we have 
$H_a(x) \geq H_a(a)$. Additionally, $H_a(x)$ is an unbounded function of $x$. Therefore,
\begin{align*}
    V(\BoundarySet^{(i)}) 
    &= \{V(\vx) \mid \vx \in \cXbar, x_i \in \{0,1\}\} = ( C_{\vpi}+1+\sum_{k\in[n]_{-i}} H_{h(\pi_i,\pi_k)}(h(\pi_i,\pi_k)),\infty ).
\end{align*}
Note that ${\Wcal_M = \{\vx \in \cX \mid V(\vx) \leq M\}\subset \cX=(0,1)^n}$ is not closed for any level 
\[
M > M_{\min}:=\min_{i\in[n]} C_{\vpi}+ \log 2 + \sum_{k\in [n]_{-i}} H_{h(\pi_i,\pi_k)}(h(\pi_i,\pi_k)),
\]
as for such an $M$, $V(\vx)< M$ for some $\vx\in \BoundarySet^{(i)}$ for some $i\in [n]$ but as a sublevel set in $\cX$, $\vx\not \in \Wcal_M$.
%

In the following lemma, we prove that the Lyapunov function $V(\vx)$ satisfies certain properties with the function $\vf(\vx)$.  For any $ M >0$ we denote the sublevel set of $V(\vx)$ over $\cXbar$ as ${\Wbar_M := \{\vx \in \cXbar \mid V(\vx) \leq M\}}.$

\begin{lem}\label{lem:V_asmp_satisfied}
    The Lyapunov-candidate  ${V:\cXbar \to [0,\infty)}$, \av{defined through \eqref{eqn:lyap}}, 
    for the vector-field $\vf:\cXbar \to \R^n$ satisfies the following properties:
    \begin{enumerate}[label=(\roman*).]
        \item $\inp{\nabla V(\vx)}{\vf(\vx)} \leq 0$  for any $\vx \in \cXbar$;
        \item there exists $\bar{M}_0$ such that 
        $$\Ebar := \Ecal \cup \Eboundary \subseteq \{\vx\in \cXbar \mid V(\vx) < \bar{M}_0\};$$
        \item for any $\bar{M}_1 \in (\Mbar_0,\infty)$, $\bar{\Wcal}_{\bar{M}_1}$ is a compact set; 
        \item the closure of $V(\Ebar)$ has an empty interior. 
    \end{enumerate}
\end{lem}
The proof of Lemma~\ref{lem:V_asmp_satisfied} is stated in Appendix~\ref{appendix:proof_extreme}.
\subsection{Convergence of Estimator}
In the following lemma, we show that the total error accumulated 
through the zero-difference martingale term converges almost surely. 

\begin{lem}\label{lem:error_convergence}
     $\sum_{t=0}^\infty\step{t} \vM(t+1)$ converges almost surely for the martingale difference sequence $\{\vM(t)\}$ defined in \eqref{eq:Ppr_def_Mt}.
\end{lem}
The proof of Lemma~\ref{lem:error_convergence} is provided in Appendix~\ref{appendix:recurrence}. We will refer to the sample paths for which $ {\sum_{t=0}^\infty \step{t}\vM({t+1}) }$ converges as $\Omegaconv$, i.e.,
\begin{align}\label{eq:omegaConv}
    \Omegaconv = \left\{\omega \in \Omega \middle| \sum_{t=0}^\infty \step{t}\vM({t+1})  \text{ converges}\right\}.
\end{align}
From Lemma~\ref{lem:error_convergence} we have $\Pr(\Omegaconv) = 1$. 

Using Lemma~\ref{lem:error_convergence}, we establish that the estimates lie in a truncation set.
\begin{lem}\label{lem:Kq_recurrence}
     For the update rule defined through \eqref{eq:projected_update_rule}, for every sample path $\omega \in \Omegaconv$, 
     there exists an index $\qrec(\omega)$ such that $\{\vPpr(t;\omega)\}$ lies in $\Kcal_{q(\omega)}$. Moreover, $\Ebar \subseteq \Kcal_{q(\omega)}$.
\end{lem}
The proof of Lemma~\ref{lem:Kq_recurrence} is provided in Appendix~\ref{appendix:recurrence}. Lemma~\ref{lem:Kq_recurrence} establishes that with probability one, $ \{\vcorrection{t}\}$ is non-zero finitely often, implying that resetting takes place finitely often.
\avnew{ To prove Lemma~\ref{lem:Kq_recurrence} we make use of Theorem~\ref{thm:2_2} which states that under certain conditions, if the Lyapunov function is below a certain threshold, then the updates will stay within a sublevel set of the Lyapunov function. The proof and the statement of Theorem~\ref{thm:2_2} closely follow from the result in \cite[Theorem 2.2]{andrieu2005stability}. However, we need to verify the proof since our results are stated for the mean-field function and the Lyapunov function being defined on an extension of the open set. The domain of the two functions in the original statement in \cite{andrieu2005stability} is an open set, which is not the case in our statements. We use the specific functions' definitions to ensure that all the sets in the proofs are well-defined and follow the desired conditions.   
The extension is essential, as otherwise the crucial property of the sublevel set of the Lyapunov function being compact is not satisfied. Note that all the proofs in this work that closely follow the proofs in \cite{andrieu2005stability} is provided in Appendix~\ref{appendix:proof_sa_main}.}
\begin{thm}\label{thm:stability_sa_main}
Consider the compact subset $\Kcal_{q(\omega)}$  of $\cXbar$ (from Lemma~\ref{lem:Kq_recurrence}). 
Then for the update rule~\eqref{eq:projected_update_rule}, we have $\limsup_{t \to \infty} d(\vPpr(t;\omega),\Ebar) = 0$ for all $\omega \in \Omegaconv$.
\end{thm}
The proof of Theorem~\ref{thm:stability_sa_main} is provided in Appendix~\ref{appendix:proof_sa_main}.

Finally, using the fact that the estimates of the update rule~\eqref{eq:projected_update_rule} lie in a truncation set and the cumulative error term convergence almost surely, we establish the almost sure convergence result of the estimates. 

\begin{pf*}{Proof of Theorem~\ref{thm:discrete_time_convergence}} 
Through Lemma~\ref{lem:error_convergence} we know that the error term $\sum_{t=0}^\infty \step{t} \vM(t+1)$ converges a.s. Also, from Lemma~\ref{lem:Kq_recurrence}, for all sample paths in the set $\omega \in \Omegaconv$
the sequence $\{\vPpr(t;\omega)\}$ stays in some compact subset $\Kcal_{\qrec(\omega)}$  of $\cXbar$. 
Moreover, Lemma~\ref{lem:Kq_recurrence} implies that $\Ebar \subseteq \Kcal_{\qrec(\omega)}$, which leads to  $\Kcal_{\qrec(\omega)}\cap \Ebar = \Ebar$. 
So, applying Theorem~\ref{thm:stability_sa_main} to every sample path in $\Omegaconv$, we conclude that $\lim_{t\to\infty}d(\vPpr(t),\Ebar)=0$ a.s.
\end{pf*}

\section{Conclusion and future work}
We presented a model for fact-checking of binary facts involving agents modeled as memoryless binary symmetric channels and proposed an online algorithm to estimate the unreliability parameters of the agents. We proved that the estimates form a dynamic process which is a stochastic approximation scheme and using results from stochastic approximation theory, we showed that it converges almost surely to the set of equilibrium points of the mean-field ODE over an extended domain $\cXbar$. In proving the convergence we studied the properties of the KL divergence used as the Lyapunov function $V(\vx)$ for the mean field ODE $\vf(\vx)$. 
The online estimator proposed in this paper and its analysis open up a variety of avenues for future work. We conjecture that the set to which the online estimator converges can be further reduced to a smaller set containing the stable equilibrium points such as $\vpi$ and $\one-\vpi$ when we have $\pi_i \neq \frac{1}{2}$ for any $i \in [n]$. Further work involves studying the convergence of variants of the proposed online estimator when not all the agents participate in the fact-checking task at every time $t$.  
\bibliographystyle{abbrv}
\bibliography{bib_new}

\begin{thebibliography}{10}

\bibitem{acemoglu2010spread}
D.~Acemoglu, A.~Ozdaglar, and A.~ParandehGheibi.
\newblock Spread of (mis) information in social networks.
\newblock {\em Games and Economic Behavior}, 70(2):194--227, 2010.

\bibitem{acemoglu2021misinformation}
D.~Acemoglu, A.~Ozdaglar, and J.~Siderius.
\newblock Misinformation: Strategic sharing, homophily, and endogenous echo chambers.
\newblock Technical report, National Bureau of Economic Research, 2021.

\bibitem{albert2004cautionary}
P.~S. Albert and L.~E. Dodd.
\newblock A cautionary note on the robustness of latent class models for estimating diagnostic error without a gold standard.
\newblock {\em Biometrics}, 60(2):427--435, 2004.

\bibitem{andrieu2005stability}
C.~Andrieu, {\'E}.~Moulines, and P.~Priouret.
\newblock Stability of stochastic approximation under verifiable conditions.
\newblock {\em SIAM Journal on control and optimization}, 44(1):283--312, 2005.

\bibitem{bonald2017minimax}
T.~Bonald and R.~Combes.
\newblock A minimax optimal algorithm for crowdsourcing.
\newblock {\em NeurIPS}, 30, 2017.

\bibitem{borkar2009stochastic}
V.~S. Borkar.
\newblock {\em Stochastic approximation: a dynamical systems viewpoint}, volume~48.
\newblock Springer, 2009.

\bibitem{budak2011limiting}
C.~Budak, D.~Agrawal, and A.~El~Abbadi.
\newblock Limiting the spread of misinformation in social networks.
\newblock In {\em Proceedings of the 20th international conference on World wide web}, pages 665--674, 2011.

\bibitem{chen2005stochastic}
H.-F. Chen.
\newblock {\em Stochastic approximation and its applications}, volume~64.
\newblock Springer Science \& Business Media, 2005.

\bibitem{chen2019beyond}
S.~Chen and A.~Moitra.
\newblock Beyond the low-degree algorithm: mixtures of subcubes and their applications.
\newblock In {\em 51st Annual ACM SIGACT Symposium on Theory of Computing}, pages 869--880, 2019.

\bibitem{dawid1979maximum}
A.~P. Dawid and A.~M. Skene.
\newblock Maximum likelihood estimation of observer error-rates using the em algorithm.
\newblock {\em J. R. Stat. Soc.: Series C (Applied Statistics)}, 28(1):20--28, 1979.

\bibitem{durrett2019probability}
R.~Durrett.
\newblock {\em Probability: theory and examples}, volume~49.
\newblock Cambridge university press, 2019.

\bibitem{feldman2008learning}
J.~Feldman, R.~O'Donnell, and R.~A. Servedio.
\newblock Learning mixtures of product distributions over discrete domains.
\newblock {\em SIAM Journal on Computing}, 37(5):1536--1564, 2008.

\bibitem{freund1999estimating}
Y.~Freund and Y.~Mansour.
\newblock Estimating a mixture of two product distributions.
\newblock In {\em Proceedings of the twelfth annual conference on Computational learning theory}, pages 53--62, 1999.

\bibitem{gao2016exact}
C.~Gao, Y.~Lu, and D.~Zhou.
\newblock Exact exponent in optimal rates for crowdsourcing.
\newblock In {\em International Conference on Machine Learning}, pages 603--611. PMLR, 2016.

\bibitem{gordon2021source}
S.~Gordon, B.~H. Mazaheri, Y.~Rabani, and L.~Schulman.
\newblock Source identification for mixtures of product distributions.
\newblock In {\em Conference on Learning Theory}, pages 2193--2216. PMLR, 2021.

\bibitem{gordon2022hadamard}
S.~L. Gordon and L.~J. Schulman.
\newblock Hadamard extensions and the identification of mixtures of product distributions.
\newblock {\em IEEE Transactions on Information Theory}, 2022.

\bibitem{hall2014martingale}
P.~Hall and C.~C. Heyde.
\newblock {\em Martingale limit theory and its application}.
\newblock Academic press, 2014.

\bibitem{hanselowski2017framework}
A.~Hanselowski and I.~Gurevych.
\newblock A framework for automated fact-checking for real-time validation of emerging claims on the web.
\newblock In {\em NIPS 2017 Workshop on Prioritising Online Content. Long Beach, USA.}, 2017.

\bibitem{hanselowski2019richly}
A.~Hanselowski, C.~Stab, C.~Schulz, Z.~Li, and I.~Gurevych.
\newblock A richly annotated corpus for different tasks in automated fact-checking.
\newblock {\em arXiv preprint arXiv:1911.01214}, 2019.

\bibitem{hui1980estimating}
S.~L. Hui and S.~D. Walter.
\newblock Estimating the error rates of diagnostic tests.
\newblock {\em Biometrics}, pages 167--171, 1980.

\bibitem{kamath2015learning}
S.~Kamath, A.~Orlitsky, D.~Pichapati, and A.~T. Suresh.
\newblock On learning distributions from their samples.
\newblock In {\em Conference on Learning Theory}, pages 1066--1100. PMLR, 2015.

\bibitem{kushner2006stochastic}
H.~Kushner and G.~Yin.
\newblock {\em Stochastic Approximation and Recursive Algorithms and Applications}.
\newblock Stochastic Modelling and Applied Probability. Springer New York, 2006.

\bibitem{nguyen2012containment}
N.~P. Nguyen, G.~Yan, M.~T. Thai, and S.~Eidenbenz.
\newblock Containment of misinformation spread in online social networks.
\newblock In {\em Proceedings of the 4th Annual ACM Web Science Conference}, pages 213--222, 2012.

\bibitem{nitzan1981characterization}
S.~Nitzan and J.~Paroush.
\newblock The characterization of decisive weighted majority rules.
\newblock {\em Economics Letters}, 7(2):119--124, 1981.

\bibitem{papanastasiou2020fake}
Y.~Papanastasiou.
\newblock Fake news propagation and detection: A sequential model.
\newblock {\em Management Science}, 66(5):1826--1846, 2020.

\bibitem{polyanskiy2014lecture}
Y.~Polyanskiy and Y.~Wu.
\newblock Lecture notes on information theory.
\newblock {\em Lecture Notes for ECE563 (UIUC)}, 6(2012-2016):7, 2014.

\bibitem{raykar2010learning}
V.~C. Raykar, S.~Yu, L.~H. Zhao, G.~H. Valadez, C.~Florin, L.~Bogoni, and L.~Moy.
\newblock Learning from crowds.
\newblock {\em Journal of machine learning research}, 11(4), 2010.

\bibitem{shapley1984optimizing}
L.~Shapley and B.~Grofman.
\newblock Optimizing group judgmental accuracy in the presence of interdependencies.
\newblock {\em Public Choice}, 43(3):329--343, 1984.

\bibitem{smyth1994inferring}
P.~Smyth, U.~Fayyad, M.~Burl, P.~Perona, and P.~Baldi.
\newblock Inferring ground truth from subjective labelling of venus images.
\newblock {\em Advances in neural information processing systems}, 7, 1994.

\bibitem{thorne2018fever}
J.~Thorne, A.~Vlachos, C.~Christodoulopoulos, and A.~Mittal.
\newblock Fever: a large-scale dataset for fact extraction and verification.
\newblock {\em arXiv preprint arXiv:1803.05355}, 2018.

\bibitem{verma2024acc}
A.~Verma, S.~Mohajer, and B.~Touri.
\newblock Distributed fact checking: Estimating unreliability.
\newblock In {\em ACC}, 2024.

\bibitem{verma2024cdc}
A.~Verma, S.~Mohajer, and B.~Touri.
\newblock Multi-agent fact-checker: Adaptive estimators.
\newblock In {\em CDC}, 2024.

\bibitem{verma2023binary}
A.~Verma, A.~Sharbafchi, B.~Touri, and S.~Mohajer.
\newblock Distributed fact checking.
\newblock In {\em ISIT}, 2023.

\bibitem{zhang2016spectral}
Y.~Zhang, X.~Chen, D.~Zhou, and M.~I. Jordan.
\newblock Spectral methods meet em: A provably optimal algorithm for crowdsourcing.
\newblock {\em The Journal of Machine Learning Research}, 17(1):3537--3580, 2016.

\end{thebibliography}

\appendix
\section{Characterization of $\Scal$}\label{appendix:Schar}
\begin{pf*}{Proof of Lemma~\ref{lem:Scharacterize}}
By the definition of $\mathcal{S}$ in~\eqref{eq:S_def}, $\vx \in \Scal$ iff ${g_{\vpi}(\vr) = g_{\vx}(\vr)}$ for every $\vr \in \{+1,-1\}^n$. Now, for a distinct pair $i,j\in[n]$, taking the sum of $g_{\vpi}(\vr)$ and $g_{\vx}(\vr)$  over the set 
${\{\vr \in \{+1,-1\}^n \mid r_i = r_j = +1\}}$, 
we get 
\[x_i x_j + (1-x_i)(1-x_j) = \pi_i \pi_j + (1-\pi_i)(1-\pi_j),\] 
or equivalently, 
\begin{align}\label{eq:xpi_system}
    \paran{\frac{1}{2}-x_i}\paran{\frac{1}{2}-x_j} = \paran{\frac{1}{2}-\pi_i}\paran{\frac{1}{2}-\pi_j}.
\end{align}
This equation holds for every $i,j\in[n]$. Note that the right-hand side of the above equation is strictly positive as $\pi_i,\pi_j\not=\frac{1}{2}$. 
 Now, consider $i,j,k\in [n]$. Multiplying this equation for $(i,j)$ and $(i,k)$, and dividing the result by the equation for $(j,k)$, we arrive at 
$\paran{\frac{1}{2}-x_i}^2 = \paran{\frac{1}{2}-\pi_i}^2$ for all $i\in [n]$, which leads to $\vx = \vpi $ or $\one - \vpi$. Therefore ${\Scal = \{\vpi,\one-\vpi}\}$. 
\end{pf*}

\section{Proof of Stochastic Approximation Lemmas}\label{appendix:stoch_approx}
\begin{pf*}{Proof of Lemma~\ref{lem:StochApprox}}
To prove the claim, first note that \[\step{t}\E[\vM(t+1)\mid \Fcal_t]=\E[\vec{P}(t+1)- \vec{P}(t)| \Fcal_t]- \step{t} \vf(\vec{P}(t)).\] Thus, we need to show ${\E[\vec{P}(t\hspace{-1pt}+\hspace{-1pt}1)\hspace{-1pt}-\hspace{-1pt} \vec{P}(t)| \Fcal_t]\hspace{-1pt}=\hspace{-1pt} \step{t} \vf(\vec{P}(t))}$ 
almost surely.
By definition, we know that 
$\step{t} \vftilde(\vR(t+1),\vP(t)) = \vP(t+1)-\vP(t)$
which results in 
\[{\vf(\vP(t))\!\! =\! \E[\vftilde(\vR(t+1),\vP(t)) \mid \Fcal_t]\!=\!\E_{\vR\sim g_{\vec{\pi}}}[\vftilde(\vR,\vP(t))]}\]
as $\{\vR(t)\}$ is i.i.d.\ with $\vR(t)\sim g_{\vec{\pi}}$.   

Since $\frac{L(t)-1}{L(t)+1} R_i(t+1) \in [-1,1]$ and $P_i(t) \in [0,1]$, we know that the value of the coordinates satisfy $\ftilde_i \in [-1,1]$. Therefore, $\|\vM(t+1)\|_{\infty} \leq 2$ a.s. 
\end{pf*}
\section{Proof regarding Lyapunov and ODE functions}\label{appendix:lyapunov_ode_prop}
\begin{pf*}{Proof of Lemma~\ref{lem:rewritten_f}}
For a fact-checker with unreliability vector $\vx \in \cX$, we define the log-odds value associated with agent $i$ as $\ell_{x_i} := \log \frac{1-x_i}{x_i}$ and the log-odds vector as $\bl_{\vx} = \transpose{(\ell_{x_1},\ell_{x_2},\dots,\ell_{x_n})}$. The probability of the output vector being $\vr$ for a set of agents with unreliability vector $\vx $, $\Pr(\vR = \vr; \vx) $ is given by
\begin{align*}
     g_{\vx}(\vec{r}) &=\frac{1}{2}\prod_{i=1}^n x_i^{\frac{1+r_i}{2}} (1-x_i)^{\frac{1-r_i}{2}} \hspace{-0.2cm} +  \frac{1}{2}\prod_{i=1}^n x_i^{\frac{1-r_i}{2}} (1-x_i)^{\frac{1+r_i}{2}}\\
     &=\frac{1}{2} \exp \paran{\sum_{i=1}^n \frac{1+r_i}{2} \log x_i +\frac{1-r_i}{2} \log (1-x_i)} +
     \frac{1}{2} \exp \paran{\sum_{i=1}^n \frac{1-r_i}{2} \log x_i +\frac{1+r_i}{2} \log (1-x_i)\hspace{-2pt}} \cr 
     &=\frac{1}{2} \exp \paran{\sum_{i=1}^n\frac{r_i}{2} \log \frac{x_i}{1-x_i} +\frac{1}{2}\log x_i(1-x_i)} 
     +\frac{1}{2} \exp \paran{\sum_{i=1}^n\frac{-r_i}{2} \log \frac{x_i}{1-x_i} +\frac{1}{2}\log x_i(1-x_i)\hspace{-2pt}} \cr 
     &= \hspace{-2pt}\paran{\hspace{-1pt}\frac{1}{2} \exp {\frac{-\inp{\vr}{\bl_{\vx}}}{2}}  \hspace{-1pt}+ \hspace{-1pt}\frac{1}{2} \exp {\frac{\inp{\vr}{\bl_{\vx}}}{2}} \hspace{-1pt}}\hspace{-1pt}
      \prod_{i=1}^n \sqrt{x_i(1\hspace{-1pt}-\hspace{-1pt}x_i)}.
\end{align*}
where we used the log-odds vector $\bl_{\vx}$ to express the probabilities in the last step.
Then, using $L = \exp \paran{-\inp{\bl_{\vx}}{\vr}} $, the function $\vftilde(\vr,\vx) $ can be rewritten as 
    \begin{align}
       \vftilde(\vr,\vx) &=  
        \frac{1}{2}
        \paran{
           \one + 
           \frac{\exp\paran{-\inp{\vr}{\bl_{\vx}}}-1}
           {\exp\paran{-\inp{\vr}{\bl_{\vx}}}+1} \vr 
        } - \vx 
        = \frac{1}{2}
        \paran{
           \one -\tanh\paran{  \frac{1}{2} \inp{\vr}{\bl_{\vx}} } \vr 
        } - \vx.\nonumber
    \end{align}
\end{pf*}
\begin{pf*}{Proof of Theorem~
\ref{thm:KL_lyapunov}}
The KL Divergence between $g_{\vec{\pi}}$ and $g_{\vx}$ can be expressed as 
\begin{align*}
    V(\vx) =  \E_{\vR\sim g_{\vec{\pi}}}[\log g_{\vec{\pi}}(\vR)] 
    - \E_{\vR\sim g_{\vec{\pi}}}[\log g_{\vx}(\vR)].
\end{align*}
Therefore for any $i\in [n]$ the partial derivative of $V(\vx)$ with respect to $x_i$  is given by 
$$\del{V(\vx)}{x_i}\! =\!\! -\del{\E_{\vR\sim g_{\vec{\pi}}}[\log g_{\vx}(\vR)] }{x_i}
\!=\!\! -\E_{\vR\sim g_{\vec{\pi}}} \left[\del{\log g_{\vx}(\vR) }{x_i} \right]\!\!,$$
where the second equality follows due to the finite support of random variable $\vR$. 

From \eqref{eq:output_prob_alternate} we know that 
\begin{align}
\del{\log g_{\vx}(\vr)}{x_i}
&=\del{}{x_i}\left\{\frac{1}{2}\sum_{i=1}^n \log(x_i(1-x_i)) + \log \paran{ \frac{1}{2} \exp \frac{-\inp{\vr}{\bl_{\vx}}}{2} + \frac{1}{2} \exp \frac{\inp{\vr}{\bl_{\vx}}}{2}}\right\}\cr
 &= \frac{1}{2 x_i(1-x_i)}\paran{1-2x_i - r_i\frac{ \exp {\frac{\inp{\vr}{\bl_{\vx}}}{2}}  -\exp {\frac{-\inp{\vr}{\bl_{\vx}}}{2}}}{ \exp {\frac{\inp{\vr}{\bl_{\vx}}}{2}}  +\exp {\frac{-\inp{\vr}{\bl_{\vx}}}{2}}} }\cr 
&= \frac{1}{x_i(1-x_i)} \paran{
\frac{1}{2} \paran{1-r_i\tanh\frac{\inp{\vr}{\bl_{\vx}}}{2}} - x_i
} \nonumber\\
&= \frac{1}{x_i(1-x_i)} \ftilde_i(\vr,\vx). \nonumber
\end{align}
Also, for any $i\in [n]$,
\begin{align}
    \del{V(\vx)}{x_i} &=  -\E_{\vR\sim g_{\vec{\pi}}} \left[\del{\log g_{\vec{x}}(\vx)}{x_i} \right]= -\E  \left[\frac{1}{x_i(1-x_i)} \ftilde_i(\vR,\vx) \right]= -\frac{1}{x_i(1-x_i)} f_i(\vx). \nonumber
\end{align}
For $\vx\in \cX$ the derivative of the trajectory $\dot \vx = \vf(\vx)$ is 
\begin{align}
    &\inp{\nabla V(\vx)}{\vf(\vx)} = -\sum_{i=1}^n \frac{1}{x_i(1-x_i)}(f_i(\vx))^2 =-\sum_{i=1}^n {x_i(1-x_i)}\paran{\del{V(\vx)}{x_i}}^2 \leq 0, \label{eq:dVdt_delV}
\end{align}
where equality holds if and only if $f_i(\vx) =0$ for all $i\in [n]$. Therefore we have 
$$\{\vx \in \cX \mid \inp{\nabla V(\vx)}{\vf(\vx)} =0  \} =\{\vx \in \cX \mid \vf(\vx) = \vec{0} \}. $$
\end{pf*}
\vspace{-10pt}
\section{Proof of Results on Extreme Behavior}\label{appendix:proof_extreme}
\vspace{-10pt}
\begin{pf*}{Proof of Lemma~\ref{lem:Vvalue_extremepoint}} 
    The Lyapunov function $V(\vx)$ is given by
    $${V(\vx) = \E_{\vR \sim g_{\vpi}} [\log g_{\vpi}(\vR)] -  \E_{\vR \sim g_{\vpi}} [\log g_{\vx}(\vR)]}.$$ 
    For brevity, we drop the distribution from the expectation notation, \soh{and simply write} $\E_{\vR}[\cdot]$.
    
    Without loss of generality assume that \soh{$\vx \in \BoundarySet^{(1)}$ and ${x_1=0}$,} which imply  ${x_i \in (0,1) }$ for all $i\in [n]_{-1}$. Then, the probability of the output vector $\vR$ is given by 
    \begin{align*}
        g_{\vx}(\vr) &= \soh{\Pr(\vR =\vr ; \vx) }\cr
        &=\soh{\sum_{s\in\{-1,+1\}}  \Pr(S=R_1=r_1, \{R_i = r_i, \forall i\in [n]_{-1}\})} \cr 
        &=\hspace{-1pt} \frac{(1\!+\!r_1)}{4} \prod_{i=2}^n x_i^{\frac{1\!-\!r_i}{2}}(1-x_i)^{\frac{1\!+\!r_i}{2}}  +  \frac{(1-r_1)}{4} \prod_{i=2}^n x_i^{\frac{1\!+\!r_i}{2}}(1-x_i)^{\frac{1-r_i}{2}}.
    \end{align*}
    By \cite[Theorem 4.1.13]{durrett2019probability} (the tower rule), we get
    \begin{align*}
        \E_{\vR}[\log g_{\vx}(\vR)]  &= 
        \E_{R_1}[\E_{\vR\mid R_1}[\log g_{\vx}(\vR)]\mid R_1] 
        = \E_{\vR\mid R_1=+1}[\log g_{\vx}(\vR)]\mid R_1 = +1],
    \end{align*}
    where the last equality follows due to $g_{\vx}(\vr)= g_{\vx}(-\vr)$ and $\Pr(R_1 = +1) = \frac{1}{2}$.  
    Therefore, 
    \soh{
    \begin{align*}
        \E_{\vR|R_1=+1}[\log g_{\vx}(\vR)]  
        &=\! - \! \log 2 \! + \! \sum_{i=2}^n \E_{\vR|R_1=+1} \!\!\left[\hspace{-1pt}{\frac{1\!-\!R_i}{2}\log x_i \hspace{-2pt}+\hspace{-2pt} \frac{1\!+\!R_i}{2}\log (1\hspace{-2pt}-\hspace{-2pt}x_i)}\hspace{-1pt}\right] \cr 
        &=\! - \! \log 2 \! + \! \sum_{i=2}^n  \left[\frac{1-\E_{\vR|R_1=+1}[R_i]}{2}\log x_i \right.\cr
        & \hspace{70pt} + \left.\frac{1+\E_{\vR|R_1=+1}[R_i]}{2}\log (1-x_i)\right] \cr 
        &= -\log 2+
        \sum_{i=2}^n  h(\pi_1,\bpi_i)\log{x_i} + h(\pi_1,\pi_i)\log (1-x_i),
    \end{align*}
    where the last equality follows from 
    $$\E_{\vR|R_1=+1}[R_i] = h(\pi_1,\pi_i) - h(\pi_1,\bar{\pi}_i)= 2h(\pi_1,\pi_i)-1 = 1- 2 h(\pi_1,\bar{\pi}_i).$$
    }
    Finally, due to the symmetry $g_{\vx}(\vr) = g_{1-\vx}(\vr)$, \soh{we get}
    \begin{align*}
        &V(\vx)  \soh{= \E_{\vR \sim g_{\vpi}} [\log g_{\vpi}(\vR)] -  \E_{\vR \sim g_{\vpi}} [\log g_{1-\vx}(\vR)]}\cr
        &
        = C_{\vpi}\!+\!\log 2 
         -\!\!\sum_{i=2}^n (h(\pi_1,\pi_i)\log x_i  + h(\pi_1,\bpi_i)\log{(1-x_i)})\cr
        &= \soh{C_{\vpi} +\log 2} 
        -  \sum_{i=2}^n (h(\pi_1,\pi_i)\log h(x_1,x_i)  + h(\pi_1,\bpi_i)\log{ (1-h(x_1,x_i)})
    \end{align*}
    for any $\vx \in \BoundarySet$ with $x_1=1$. Compactly expressed for $\vx \in \BoundarySet^{(1)}$ we have
    \begin{align}\label{eq:Vdef_extreme}
        \soh{V(\vx) = C_{\vpi}+\log 2 + \sum_{k=2}^n H_{h(\pi_1,\pi_k)}h(x_1,x_k),}
    \end{align} 
    \soh{which proves the first claim of the lemma.}
    
\soh{Next, for $\vx \in \BoundarySet$ with $x_1=0$, from~\eqref{eqn:LratioExtension} we have
$\frac{L(\vR,\vx)-1}{L(\vR,\vx)+1}=-R_1$, and $\vf(\vx) = \frac{1}{2}\E_{\vR\sim g_{\vpi}}[\one-R_1 \vR] - \vx$.  This yields $\f_1(\vx)\hspace{-1pt}=\hspace{-1pt} 0$. Also, for any $i\hspace{-1pt}\in\hspace{-1pt} [n]_{-1}$,  we get
    \begin{align*}
    \f_i(\vx)\!=\! \frac{1}{2}\E[1\!\!-\!\!R_1R_i]\hspace{-1pt} \!\!-\!\!x_i= \Pr(R_1\neq R_i)\!\!-\!\!x_i = \hfunc(\pi_1,\bar{\pi}_i)\!\!-\!\!x_i.    
    \end{align*}
    Similarly, if $\vx \in \BoundarySet$ and $x_1 = 1$, we have ${\frac{L(\vR,\vx)-1}{L(\vR,\vx)+1}= R_1}$ and ${\vf(\vx) = \frac{1}{2}\E_{\vR\sim g_{\vpi}}[\one+R_1 \vR] - \vx }$. Therefore, ${\f_1(\vx)= 0}$ and ${\f_i(\vx)\!=\! \frac{1}{2}\E[1+R_1R_i]\hspace{-1pt} -x_i=  \hfunc(\pi_1,{\pi}_i)-x_i}$. These complete the proof of the second claim. 
}

\soh{Lastly, for $\vx \in \BoundarySet^{(1)}$ and $i\in [n]_{-1}$, from part~(i) of the lemma we have 
    \begin{align*}
        \del{V(\vx)}{x_i} 
        &=-\frac{\hfunc(\pi_1,\hfunc(x_1,\pi_i))}{x_i} + \frac{\hfunc(\pi_1,1-\hfunc(x_1,\pi_i))}{1-x_i} 
        = \frac{x_i-\hfunc(\pi_1,\hfunc(x_1,\pi_i))}{x_i(1-x_i)}.
    \end{align*}
}
    On the other hand for $i=1$, since $\del{V(\vx)}{x_1}$ is finite, we have $f_1(\vx) =-x_1(1-x_1)\del{V(\vx)}{x_1} = 0$ for $\vx \in \BoundarySet^{(1)}$. 

  Hence,  for every $\vx \in \BoundarySet$ and every $i\in [n] $ we have 
    ${f_i(\vx) \hspace{-1pt}=\hspace{-1pt} -x_i(1\hspace{-1pt}
-\hspace{-1pt}x_i)\del{V(\vx)}{x_i}}$.
\end{pf*}
%
\vspace{-5pt}
\begin{pf*}{Proof of Lemma~\ref{lem:V_asmp_satisfied}}
    From Theorem~\ref{thm:KL_lyapunov}, for any $\vx \in \cX$, we know that $\inp{\nabla V(\vx)}{\vf(\vx)}\leq 0 $. Moreover,
     from Lemma~\ref{lem:Vvalue_extremepoint}, for  $\vx \in \BoundarySet $, we have $f_i(\vx) = x_i(1-x_i)\del{V(\vx)}{x_i}$ for all $i\in [n]$. Therefore for $\vx \in \BoundarySet$, \eqref{eq:dVdt_delV} holds and $\inp{\nabla V(\vx)}{\vf(\vx)} \leq 0$. Thus we have property $(i)$. 

 To prove (ii) consider $\vx^{\ast}\in \cXbar$ to be any zero of the function $\vf(\vx)$. For any two distinct $i,j \in [n]$, we have $\frac{f_i(\vx^{\ast})+f_j(\vx^{\ast})}{2}=0$. This leads to 
\begin{align}
    1-x_i^{\ast}-x_j^{\ast} &= \E\left[\tanh\left(\frac{1}{2}\inp{\vR}{\ell_{\vx^{\ast}}}\right)\frac{(R_i+R_j)}{2}\right] \cr
    &\leq \E\left[\left|\frac{R_i+R_j}{2}\right|\right] = \E[\indicator{R_i=R_j}] = h(\pi_i,\pi_j), \nonumber
\end{align}
where the inequality holds because $|\tanh a|\leq 1$ for any ${a\in \R}$. Rearranging the inequality, we get ${x_i^{\ast}+x_j^{\ast} \geq h(\pi_i,\bpi_j)}$.

Similarly, considering 
$\left|\frac{f_i(\vx^{\ast})-f_j(\vx^{\ast})}{2}\right|=0$ for any distinct $i,j\in[n]$,
we obtain  ${|x_i^{\ast}-x_j^{\ast}| \leq h(\pi_i,\bpi_j).}$

In summary, for any zero $\vx^{\ast}$ of the function $\vf(\vx)$, for any $i, j \in [n]$, we have 
$|x_i-x_j| \leq h(\pi_i,\bpi_j) \leq x_i+x_j$. 
Thus,
\begin{align}\label{ineq:equil_subset}
    \{\vx \in \cXbar \mid \vf(\vx) = \vec{0} \} \subseteq \mathcal{A},    
\end{align}
where the set $\mathcal{A}$ is defined as 
\begin{align*}
\mathcal{A} \hspace{-1pt}:=\hspace{-1pt}\{\vx\in\cXbar \mid |x_i\hspace{-1pt}-\hspace{-1pt}x_j| \leq h(\pi_i,\bpi_j) \hspace{-1pt}\leq\hspace{-1pt} x_i\hspace{-1pt}+\hspace{-1pt}x_j, \forall i,j \in [n]\}.
\end{align*}
Note that $\mathcal{A}$ is compact and a proper subset of $\cXbar$. From the definition of KL divergence, we know that 
$
{\sup_{\vx\in\mathcal{A}}V(\vx) < \Mbar_0.}
$
for a finite constant $\Mbar_0.$ Therefore,
$
\mathcal{A}\subseteq \{\vx \in \cX \mid V(\vx) < \Mbar_0\},
$
which along with \eqref{ineq:equil_subset} gives us $\{\vx\in \cXbar \mid \vf(\vx) = \vec{0}\}\subseteq \{\vx\in\cX \mid V(\vx)< \Mbar_0\}$.
%
    Therefore, 
    ${\Ebar \subseteq \{\vx \in \cXbar \mid V(\vx) < \Mbar_0\}.}$ 

    For any $C >0 $, we know that $\Wbar_{C}$ is a closed subset of $\cXbar$.
    Therefore, $\Wbar_{C}$ is a compact set.  
    
    Finally, from Sard's theorem, we know that the closure of $V(\Ecal)$ has an empty interior since 
    ${\Ecal = \{\vx \in \cX \mid \nabla V(\vx) = 0\}}$. 
    On the other hand, $\Eboundary$ has finitely many isolated points. Thus, the closure of $V(\Ecal \cup \Eboundary)$ has an empty interior. 
    %
\end{pf*}
\vspace{-5pt}
\section{Proof of Recurrence}\label{appendix:recurrence}
\vspace{-5pt}
\begin{pf*}{Proof of Lemma~\ref{lem:error_convergence}}
    For any $\vr \hspace{-1pt}\in \hspace{-1pt}\{\pm 1\}^n$ and $\vx \hspace{-1pt}\in \hspace{-1pt}\cX$ we know that 
    $\|\vftilde(\vr,\vx) - \vf(\vx)\|^2 \leq 4n$. So we have 
    $\E[\| \step{t}\vM({t+1})\|^2 \mid \Fcal_t] \leq 4n\step{t}^2$.

    Since $\step{t}\vM(t+1)$ is a zero-difference martingale, the sum of the co-ordinates of $\step{t}\vM(t+1)$, \emph{i.e.}, $\sum_{\ell=0}^t \step{\ell}M_i(\ell+1)$, also form a martingale.
    Since $\sum_{t=0}^\infty \step{t}^2 <\infty$, by the convergence theorem in \cite[Theorem 2.17]{hall2014martingale} for any coordinate $i\in [n]$, we know that $\sum_{t=0}^\infty \step{t}M_i(t+1)$ converges a.s. Therefore, $  {\sum_{t=0}^\infty \step{t}\vM(t+1)}$ converges a.s.
\end{pf*}

In order to prove the recurrence of estimates in a compact set, we need the following result. The following theorem ensures that after a large enough time such that the step-sizes and the accumulated error from that time forward are small enough, then the estimates stay within a sublevel set of the Lyapunov function. The following theorem is based on Theorem 2.2 in {\cite{andrieu2005stability}}. 
\begin{thm}\label{thm:2_2}
Consider the function $\vf:\cXbar \to \cXbar$ and $V:\cXbar \to [0,\infty)$ defined in \eqref{eq:fi_def} and \eqref{eqn:lyap}, respectively.
For any ${\Mbar \in (\Mbar_0,\Mbar_1]}$,  where $\Mbar_0$ is the constant in Lemma~\ref{lem:V_asmp_satisfied} and $\Mbar_1 \in (\Mbar_0,\infty)$, there exist $\delta_0$, $\lambda_0 \in \R^+$ such that for all $t\geq 1$, all ${\vtheta(0) \in \Wbar_{\Mbar_0}}$, all sequences $\{\step{t}\}$ of non-negative numbers, and all sequences $\{\vdeterror(t)\}$ in $[-1,1]^n$
satisfying 
\begin{align*}
    &\sup_{0\leq k \leq t} \step{k} \leq \lambda_0, 
    \quad 
    \sup_{0\leq k \leq t} \norm{ \sum_{\ell=0}^k \step{\ell} \vdeterror(\ell+1) } \leq \delta_0,
\end{align*}
and $\vtheta(k) \in \cXbar$  for all $k\in [t]$,
with
$\vtheta(k) = \vtheta(k-1)  + \step{k-1}(\vf(\vtheta(k-1)) + \vdeterror(k))$, 
we have $V(\vtheta(k))\leq \Mbar$  for $k \in [t]$.
\end{thm}
Using the above result, we prove Lemma~\ref{lem:Kq_recurrence}.
\begin{pf*}{Proof of Lemma~\ref{lem:Kq_recurrence}}
We prove the result for every sample path 
$${\omega \in \Omegaconv = \{\omega \in \Omega \mid \sum_{t=0}^{\infty} \step{t}\vM(t+1;\omega) \text{ converges}\}}.$$
For brevity, we drop the $\omega$ from the notation of random variable unless needed for clarity.  

Define the supremum of the function $V(\cdot)$ over the initial set $\Kcal_0$ as $C_0 := \sup \{V(\vx) \mid \vx \in \Kcal_0 \}$.
We know that $C_0 < \infty$ since $\Kcal_0$ is a compact subset of $\cXbar$ and $V(\vx)<\infty$ for all $\vx \in \cXbar.$ Also, let $C_M := \max \{C_0, \Mbar_0\}$, where $\Mbar_0$ is the constant in Lemma~\ref{lem:V_asmp_satisfied}.
$C_M$ represents the constant for which the sublevel set of $V$ contains the union of the equilibrium set $\Ebar$ and the initial truncation set $\Kcal_0$. Recall that for the resetting estimate we have $\vPinit \in \Kcal_0$. Therefore, we have $\Ebar \cup \Kcal_0 \subseteq \Wbar_{C_M}$.  
Recall that for any ${C>0}$, $\Wbar_{C}$ is a compact subset of $\cXbar$ and $\cup_{t=0}^{\infty} \Kcal_t = \cXbar$.
So we know for any ${C_M' \in (C_M,\infty)}$ there exists $\qrec \in \N_0$ such that 
${\Wbar_{C_M'} \subseteq \Kcal_{\qrec} }$. 

Fix $t_0 \geq 0$. Then, applying Theorem~\ref{thm:2_2} to the sequence started at time $t_0$, there exists $\delta_0, \lambda_0 \in \R^+$ such that for all sequences $\{\vPpr(t)\}$ with $\vPpr(t_0) \in \Wcal_{C_M}$, if we have 
\begin{align}\label{eq:conditions868}
    \sup_{t_0 \leq k \leq t} \step{k} \leq \lambda_0,  
    \sup_{t_{0} \leq k \leq t} \norm{ \sum_{\ell=t_0}^k \step{\ell}\vM(\ell+1) } \leq \delta_0,
\end{align}
for some $t \geq t_0 + 1$, then  $V(\vPpr(k)) \leq C_M'$ for ${k \in \{t_0,\dots,t\}}$. 

Contrary to the conclusion, let us assume that resetting occurs infinitely often.

For $t \geq t_0$, define $\Treset(t) := \min\{k \leq t : \proj{k} =\proj{t} \} $ as the most recent index at which resetting takes place. 

We have $\lim_{t\to \infty }\|\sum_{k=t}^\infty \step{k}\vM(k+1)\| =0$ a.s. since from Lemma~\ref{lem:error_convergence} we know that $\sum_{k=0}^\infty \step{k}\vM(k+1)$ converges a.s.
Therefore, there exists a $T_0$ after which the step-sizes sequence $ \{\step{t} \mid t \geq T_0\}$ and the zero-difference martingale sequence $\{\vM(t+1) \mid t\geq T_0\}$ satisfy $\eqref{eq:conditions868}$.

Since resettings occur infinitely often, we know there exists $T_1$ such that 
$ \Treset(T_1-1) > T_0 $, ${T_1 = \Treset(T_1)}$, and $\proj{T_1} > q$. 
In other words, $T_1$ is a time after $T_0$ by which two resettings have taken place. Define 
\begin{align}\label{eq:yvecProof7}
    \vec{y} = \vPpr(T_1-1)+ \step{T_1-1}\vftilde(\vR(T_1),\vPpr(T_1-1)).
\end{align}
We know $\vec{y}  \not \in \Kcal_{\qrec}$  since
$\proj{T_1} \geq \qrec+1$ and $\Treset(T_1) = T_1$. Furthermore, ${\vPpr(\Treset(T_1\hspace{-1pt}-\hspace{-1pt}1)) \hspace{-1pt}= \hspace{-1pt}\vPinit \hspace{-1pt}\in\hspace{-1pt} \Kcal_0}$ for the resetting time $\Treset(T_1\hspace{-1pt}-\hspace{-1pt}1)$. Since ${\Treset(T_1\hspace{-1pt}-\hspace{-1pt}1)\hspace{-1pt}>\hspace{-1pt}T_0}$ we get ${\{\step{t} \mid t \hspace{-1pt}\geq\hspace{-1pt} \Treset(T_1)\}}$ and ${\{\vM(t+1) \hspace{-1pt}\mid t\hspace{-1pt}\geq\hspace{-1pt} \Treset(T_1)\}}$ satisfy~\eqref{eq:conditions868}. From Theorem~\ref{thm:2_2} we know that for $\vec{y}$, defined in \eqref{eq:yvecProof7}, ${V(\vec{y}) \leq C_M'}$ which contradicts 
${\vec{y} \not\in \Kcal_{\qrec}}$ since ${\Wcal_{C_M'} \subseteq \Kcal_q}$.
\end{pf*}

%
\section{Proof of Stochastic Approximation Result }\label{appendix:proof_sa_main}
The proof in this section closely follows that of~\cite{andrieu2005stability}. We begin with showing some properties for functions $\vf$ and $V$. 

\begin{lem}\label{lem:2_1}
    Consider the functions $\vf:\cXbar \to \cXbar$ and ${V:\cXbar \to [0,\infty)}$ defined through \eqref{eq:fi_def} and \eqref{eqn:lyap} respectively (with function definitions extended to $\cXbar$).
    \vspace{-10pt}
    \begin{enumerate}[wide, labelwidth=!, labelindent=0pt, label=(\roman*)]
        \item Let $\Kcal \subset \cXbar$ be a subset such that $0 < \inf_{\vtheta \in \Kcal} |\inp{\nabla V}{\vf}| $. For any $0< \delta <  \inf_{\vtheta \in \Kcal} |\inp{\nabla V}{\vf}| $, there exist $\lambda >0$ and $\beta > 0 $ such that, for any $\rho \in [0,\lambda]$, $\vzeta $, $\norm{\vzeta} \leq \beta$, and $\vtheta \in \Kcal$, $V(\vtheta + \rho \vf(\vtheta)+ \rho\vzeta) \leq V(\vtheta) - \rho\delta$. 
        \item For any $\Mbar \in (\Mbar_0,\Mbar_1]$, where $M_0,M_1$ are defined in Lemma~\ref{lem:V_asmp_satisfied}, there exist $\lambda >0$ and $\beta > 0 $ such that, for any $\rho \in [0,\lambda]$, $\vzeta $, $\norm{\vzeta} \leq \beta$, and $\vtheta \in \Wbar_{\Mbar}$, ${\vtheta + \rho \vf(\vtheta)+ \rho\vzeta \in \Wbar_{\Mbar}}$.  
    \end{enumerate}
    \vspace{-10pt}
\end{lem}
\begin{pf}
\vspace{-10pt}
\begin{enumerate}[wide, labelwidth=!, labelindent=0pt, label=(\roman*)]
    \item  For any $0 < \delta <\inf_{\vtheta \in \Kcal} |\inp{\nabla V}{\vf}| $, there exist $\lambda > 0$ and $\beta > 0 $ such that for all $\rho \in [0,\lambda]$ and $\norm{\vzeta}\leq  \beta$, and $t \in [0,1]$ we have $\vtheta \in \Kcal$, $\vtheta + \rho t \vf(\vtheta) +\rho t \vzeta \in \cXbar$ \avproof{(existence ensured since $\vf$ will be finite over compact set $\Kcal$ since $\vf$ is bounded)}, and 
    \begin{align}
     &|\inp{\nabla V(\vtheta)}{\vf(\vtheta)} - \inp{\nabla V(\vtheta + \rho t \vf(\vtheta) +\rho t \vzeta )}{\vf(\vtheta) + \vzeta}  \cr 
      &\leq \inf_{\vtheta \in \Kcal} |\inp{\nabla V}{\vf}| - \delta.  \nonumber
    \end{align}

    We know 
    \begin{align*}
        &V(\vtheta + \rho  \vf(\vtheta) +\rho  \vzeta) - V(\vtheta) \cr 
        &= \rho \int_{0}^1  \inp{\nabla V(\vtheta + \rho t \vf(\vtheta) +\rho t \vzeta )}{\vf(\vtheta) + \vzeta} dt \cr
        &=\rho \inp{\nabla V(\vtheta)}{\vf(\vtheta)} 
        + \rho \int_{0}^1(  \inp{\nabla V(\vtheta + \rho t \vf(\vtheta) +\rho t \vzeta )}{\vf(\vtheta) + \vzeta}- \inp{\nabla V(\vtheta)}{\vf(\vtheta)} ) dt \leq - \rho\delta.
    \end{align*}
    \item Consider $\Mbar' \in (\Mbar_0,\Mbar)$. Since $\vf$ is bounded and $V$ is continuous, there exists $\lambda_0, \beta_0 \in \R^+$ such that for all $\rho \in [0,\lambda_0]$ and $\norm{\vzeta} \leq \beta_0$, and ${\vtheta \in \Wbar_{\Mbar'}}$ we have $\vtheta + \rho  \vf(\vtheta) +\rho  \vzeta \in \Wbar_{\Mbar}$.
    Applying the result of $(i)$ to the set
    \begin{align*}
    \Kcal &= \{\vtheta \in \cXbar | \Mbar' \leq V(\vtheta) \leq \Mbar \} = \Wcal_{\Mbar}\setminus \{\vtheta\in \cXbar | V(\vtheta)<\Mbar'\},    
    \end{align*}
     there exists  $\lambda_1, \beta_1 \in \R^+$ such that for all 
    ${\rho \in [0,\lambda_1]}$ and $\norm{\vzeta}\leq \beta_1$, and $\vtheta \in \Kcal$, we have 
    $$V(\vtheta + \rho  \vf(\vtheta) +\rho  \vzeta) \leq V(\vtheta) \leq \Mbar$$
    implying
    $\vtheta + \rho  \vf(\vtheta) +\rho  \vzeta \in \Wbar_{\Mbar}$. 
\end{enumerate}
\end{pf}
Now we provide the proof for Theorem~\ref{thm:2_2} following the proof of Theorem 2.2 in ~\cite{andrieu2005stability}. 
\begin{pf*}{Proof of Theorem \ref{thm:2_2}}
Consider some $\Mbar' \in (\Mbar_0,\Mbar)$. From Lemma~\ref{lem:2_1}, we know that there exists $\lambda_0,\beta_0 \in \R^+$ such that for all $\vtheta, \rho$, and $ \vzeta$ satisfying $V(\vtheta) \leq \Mbar'$, ${\rho \in [0,\lambda_0]}$, and $\norm{\vzeta}\leq \beta_0$, we have ${V(\vtheta + \rho \vf(\vtheta) + \rho \vzeta) \leq \Mbar'.}$

By continuity of $\vf$ and $V$ there exists $\delta_0\in (0,\beta_0]$ such that for all $\vtheta \times {\vtheta}' \in \cXbar \times \cXbar$ satisfying $V(\vtheta)\leq \Mbar $ and $\norm{\vtheta - \vtheta' }\leq \delta_0$, we have  
\begin{align}\label{eq:continuity1}
    \norm{ \vf(\vtheta )\hspace{-1pt} - \hspace{-1pt}\vf(\vtheta') } \hspace{-1pt}\leq\hspace{-1pt} \beta_0 \text{ and }
    |V(\vtheta)\hspace{-1pt}-\hspace{-1pt} V(\vtheta')| \hspace{-1pt}\leq\hspace{-1pt} \Mbar \hspace{-1pt}-\hspace{-1pt} \Mbar'. 
\end{align}
We will use induction to prove for all $k \in [t]$, we have $V(\vtheta'(k))\leq \Mbar'$, and $V(\vtheta(k))\leq \Mbar$, where the sequence $\{\vtheta'(k)\}$ is defined as $ \vtheta'(0) = \vtheta(0)$ and for all $k\in [t]$,  
    ${\vtheta'(k) 
    = 
    \vtheta'(k-1) + \step{k-1} \vf(\vtheta(k-1)).}$ 

Under the stated assumptions ${V(\vtheta'(0))= V(\vtheta(0)) \leq \Mbar_0}$. Since $0 \leq \step{0} \leq \lambda_0$ and 
$\norm{ \vtheta'(1) - \vtheta(1) } = \norm{ \step{0} \vdeterror(1) } \leq \delta_0$, on the one hand Lemma~\ref{lem:2_1} shows that 
$${V(\vtheta'(1) ) = V(\vtheta'(0) + \step{0}\vf(\vtheta(0))) \leq \Mbar'}$$
and ${ V(\vtheta'(1)) = V(\vtheta(0) + \step{0}\vf(\vtheta(0)) + \step{0}\vdeterror(1) \leq \Mbar },$
which proves the result for $t= 1$. 

Assuming the result holds for $k \in [t-1]$ for $t > 1$. By construction for $j \in [k]$, 
$\vtheta(j) - \vtheta'(j) = \vtheta(j-1) + \step{j-1}\vdeterror(j),$
which implies that 
$\vtheta(j) -  \vtheta'(j) = \sum_{i=1}^j \step{i-1}\vdeterror(i)$.

Under the stated assumptions {ensuring continuity} and \eqref{eq:continuity1}, for $j\in [k]$, we have 
$$ \norm{ \vtheta(j) -  \vtheta'(j) } \leq \delta_0 \text{ and } \norm{ \vf(\vtheta(j))-\vf(\vtheta'(j))} \leq \beta_0. $$ 
On the other hand, 
\begin{align*}
    &\vtheta'(k+1) = \vtheta'(k) + \step{k} \vf(\vtheta(k)) =  \vtheta'(k) + \step{k} \vf(\vtheta'(k)) +  \step{k} (\vf(\vtheta(k)) - \vf(\vtheta'(k)) ).
\end{align*}
Since $0 \leq \step{k}\leq \lambda_0$ and $V( \vtheta'(k)) \leq \Mbar'$, Lemma~\ref{lem:2_1} shows $V(\vtheta'(k+1))\leq M'$. 

Using ${\norm{ \vtheta(k+1) - \vtheta'(k+1)} \leq \delta_0}$, \eqref{eq:continuity1} implies that ${V(\vtheta(k+1))\leq \Mbar}$, concluding the proof.

\end{pf*}
%
Let ${A_{\delta} := \{ \vtheta \in \cXbar | d(\vtheta, A) \leq \delta\} }$ for any $A\subset \cXbar$ and  $\delta > 0$, and let 
${\norm{\phi}_A := \sup_{\vtheta \in A} \norm{\phi(\vtheta)}}$ for any function $\phi:\cXbar\!\!\! \to \!\R$.

The following Lemma based on {\cite[Lemma 2.4]{andrieu2005stability}} is used in the proof of Theorem~\ref{thm:stability_sa_main}. We will state the lemma for any $\omega \in \Omegaconv$ and drop the $\omega$-notation from $\qrec(\omega)$ for brevity. 
\begin{lem}\label{lemma:2_4}
    Under the assumption of Theorem~\ref{thm:stability_sa_main} let $\Ncal \subset \cXbar$ be a neighborhood of $\Ebar \cap \Kcal_{\qrec}$ which satisfies $\sup_{\vtheta \in \Kcal_{\qrec}\setminus \Ncal} \inp{\nabla V(\vtheta}{\vf(\vtheta)} < 0$. 
    There exist $\delta, \varepsilon, \lambda>0$ (depending on the sets $\Ncal$ and $\Kcal_{\qrec}$) such that for any $\delta' \in (0,\delta]$, $\lambda' \in (0,\lambda]$, and $\eta> 0$, one can find an integer $T$ and a sequence $\{\hvPpr(j) | j \geq T\}$ satisfying 
    \begin{align*}
        &\sup_{j\geq T} \norm{ \vPpr(j) - \hvPpr(j) } \leq \delta', 
        \quad\sup_{j\geq T} \step{j-1} \leq \lambda', \text{ and }\cr 
        &
        \,\sup_{j\geq T} |V(\vPpr(j))- V(\hvPpr(j)) | \leq \eta,
    \end{align*}
    and for $j \geq T+1$,
    \begin{align*}
        V(\hvPpr(j)) &\leq V(\hvPpr(j)) - \step{j-1}\varepsilon + (\eta+\step{j-1}\varepsilon) \indicator{\hvPpr(j-1) \in \Ncal}.
    \end{align*}
\end{lem}
\begin{pf}
     For legibility in the proof we set $\Kcal = \Kcal_{\qrec}$. 
     Let us choose $\delta_0>0$ such that the set of points in $\cXbar$ which are $\delta_0$ away from the set $\Kcal$ satisfy ${\Kcal_{\delta_0} \subset \Wbar_{\Mbar_2}\subset \cXbar}$, for some $\Mbar_2 \geq \Mbar $. The set $\Kcal_{\delta_0}\setminus \Ncal$  satisfies  ${\sup_{\vtheta \in \Kcal_{\delta_0}\setminus \Ncal} \inp{\nabla V}{\vf} < 0}$. 

     By Lemma~\ref{lem:2_1}, for any $\varepsilon >0$ such that 
     $\sup_{\vtheta \in \Kcal_{\delta_0}\setminus \Ncal} \inp{\nabla V}{\vf} < -\varepsilon$,
     one may choose $\lambda > 0 $ and $\beta > 0 $ small enough so that for any $\rho \in [0,\lambda]$ and $\norm{\vzeta} \leq \beta$, and $\vtheta \in \Kcal_{\delta_0}\setminus \Ncal$ we have
     \begin{align}\label{eq:Vineq}
         V(\vtheta + \rho \vf(\vtheta) + \rho \vzeta )  \leq V(\vtheta) - \rho \varepsilon. 
     \end{align}
     Note that {$\vf$ is bounded so $\norm{\vf}_{\Kcal}$ is finite}. So, using the uniform continuity of $\vf$ on $\Kcal$, for any $\eta > 0$ one may choose $\delta \in (0,\lambda \norm{\vf}_{\Kcal})$ small enough so that for all 
     ${(\vtheta, \vtheta') \in \Kcal_{\delta_0} \times \Kcal_{\delta_0}}$ satisfying $\norm{ \vtheta - \vtheta'} \leq \delta \leq \lambda \norm{ \vf }_{\Kcal}$ we have 
     \begin{align}\label{eq:continuity2}
         \norm{ \vf(\vtheta) - \vf(\vtheta') } \leq \beta 
         \text{ and }
         |V(\vtheta) - V(\vtheta')| \leq \eta.
     \end{align}
     Under the stated conditions (\avproof{regarding bounded step-size and bounded cumulative error}) for all $\delta' \in (0,\delta]$ and ${\lambda' \in (0,\lambda]}$ there exists an integer $T$ such that for any $t \geq T+1$, and $\step{t} \leq \lambda'$ and 
     $\norm{ \sum_{k=T+1}^t \step{k-1}\vM(k) } \leq \delta'$. 

     Define recursively for $j \geq T$, the sequence $\{\hvPpr(j) | j\geq T\}$ as $\hvPpr(T) := \vPpr(T)$ and for $j\geq T+1$, 
     \begin{align*}
         \hvPpr(j) = \hvPpr(j-1) + \step{j-1} \vf(\vPpr(j-1)).
     \end{align*}
     By construction for $j \geq T+1$, 
     \[{\hvPpr(j) - \vPpr(j) = \sum_{i=T+1}^{j} \step{i-1}\vdeterror(i)}\]
     which implies that $\sup_{j\geq T} \norm{ \hvPpr(j) - \vPpr(j) } \leq \delta'$. On the other hand, for $j\geq T+1$, 
     \begin{align*}
         \hvPpr(j) &= \hvPpr(j-1) + \step{j-1} \vf( \hvPpr(j-1) ) + 
         \step{j-1} ( \vf(\vPpr(j-1))  - \vf(\hvPpr(j-1) )),
     \end{align*}
     and since $\norm{ \hvPpr(j-1) - \vPpr(j-1) } \leq \delta' \leq \delta$, \eqref{eq:continuity2} shows that ${\norm{ \vf(\vPpr(j-1)) - \vf(\hvPpr(j-1)) } \leq \beta}$. 
     By \eqref{eq:Vineq} we know that whenever $\hvPpr(j-1) \in \Kcal_{\delta}\setminus \Ncal$ 
     since $\Kcal_{\delta} \subset \Wcal_{M_2}$,  
     ${V(\hvPpr(j))\leq V(\hvPpr(j-1)) - \step{j}\varepsilon}$. 

    Now \eqref{eq:continuity2} implies that $|V(\hvPpr(j)) -  V(\hvPpr(j-1))|\leq \eta$ for any $\vPpr(j-1) \in \Kcal_{\delta}$ and $|V(\vPpr(j)) - V(\hvPpr(j))|\leq \eta$ for any ${\vPpr(j) \in \Kcal}$.
\end{pf}
Finally, we need the following lemma from \cite{andrieu2005stability} for the proof of Theorem~\ref{thm:stability_sa_main}. 
\begin{lem}[{\cite[Lemma 2.5]{andrieu2005stability}}]\label{lemma:2_5}
Assume \ref{asmp:stepsizes} holds for a sequence $\{\step{t}\}$. Let $\varepsilon$ be a real constant, $n$ be an integer, and let $-\infty < a_1 < b_1 <\cdots < a_n < b_n < \infty$ be real numbers. 
Let $\{u_j\}$ be a bounded real sequence such that, for any $\eta > 0$, there exists an integer $J$ such that for all $j \geq J$, 
\begin{align*}
    u_j \leq 
    u_{j-1} - \step{j} \varepsilon + (\eta + \step{j}\varepsilon) \indicator{u_{j-1} \in A}, 
    \quad 
    A = \bigcup_{i=1}^n [a_i,b_i]. 
\end{align*}
Then the limit points of the sequence $\{u_j\}$ are included in A. 
\end{lem}
With all the required Lemmas we now proceed to the proof of Theorem~\ref{thm:stability_sa_main}. Since the statements of the proof hold for $\omega \in \Omegaconv$ we drop the notation of $\omega$ in the following proof. 
\begin{pf*}{Proof of Theorem~\ref{thm:stability_sa_main}}
    We first prove that $\lim_{j\to \infty} V(\vPpr(j))$ exists. 
    For any $\alpha >0$, define the set 
    $$[V(\Ebar \cap \Kcal_{\qrec})]_{\alpha}:= \{ x \in \R : d(x, V(\Ebar \cap \Kcal_{\qrec})) \leq \alpha \}.$$
\vspace{-10pt}
    Since $\Mbar_0 = \sup_{\vx\in \Ebar}V(\vx)<\infty$, we know that  
    \[{[V(\Ebar \cap \Kcal_{\qrec})]_{\alpha} \subseteq [-\alpha,\Mbar_0+\alpha]}.\] 
    By definition, $[V(\Ebar \cap \Kcal_{\qrec})]_{\alpha}$ must be a union of intervals of length at least $2\alpha$. Since $[-\alpha,\Mbar_0+\alpha]$ is a finite interval,  $[V(\Ebar \cap \Kcal_{\qrec})]_{\alpha} $
    must be a finite union of disjoint intervals of length at least equal to $2\alpha$.
    
    By Lemma~\ref{lemma:2_4}, there exist positive constants $\delta, \varepsilon, \lambda$ such that for any $\delta' \in (0,\delta]$, $\lambda' \in (0,\lambda]$, and $ \eta > 0 $, one may find an integer $T$ and a sequence 
    ${\{\hvPpr(j) | j \geq T\}}$ such that we have 
    $$\sup_{j \geq T} \norm{\vPpr(j) - \hvPpr(j)} \leq \delta \text{ and } 
    \,\, {\sup_{j\geq T} |V(\vPpr(j)) - V(\hvPpr(j))| \leq \eta.}$$
    Moreover for any $j \geq T+1$,
    \begin{align*}
        V(\hvPpr(j))& \leq V(\hvPpr(j-1)) - \step{j-1}\varepsilon + (\eta + \step{j-1} \varepsilon) \indicator{V(\hvPpr(j-1))\in [V(\Ebar\cap\Kcal_{\qrec})]_{\alpha}}, 
    \end{align*}  
    where we have chosen $\Ncal = V^{-1}(\text{int}([V(\Ebar \cap \Kcal_{\qrec})]_{\alpha}))$ and used 
    $\indicator{\vtheta \in \Ncal} \leq \indicator{V(\vtheta) \in [V(\Ebar \cap \Kcal_{\qrec})]_{\Ncal}}.$ 

    By Lemma~\ref{lemma:2_5}, the limit points of the sequence $\{V(\vPpr(j))\}_{j\geq 0}$ are in $[V(\Ebar \cap \Kcal_{\qrec})]_{\alpha'}$ for $\alpha' = \alpha + \eta$. 
    Since $\alpha$ and $\eta$ can be chosen arbitrarily small, this implies that the limit points of the sequence $\{V(\vPpr(j))\}_{j\geq 0}$ are included in ${\cap_{\alpha > 0} [V(\Ebar \cap \Kcal_{\qrec})]_{\alpha}}$. 
    We have 
    \[{V(\Ebar \cap \Kcal_{\qrec}) = \cap_{\alpha > 0} [V(\Ebar\cap\Kcal_{\qrec})]_{\alpha}}.\] 
    Thus, the limit points 
    $\{V(\vPpr(j))\}$  
    belong to the set ${V(\Ebar\cap\Kcal_{\qrec})}$. 

    But, 
    ${\limsup_{j\to\infty} |V(\vPpr(j) - V(\vPpr(j-1))| =0}$ 
    which implies that the set of limit points of $\{V(\vPpr(j))\}$ is an interval. 
    Because $V(\Ebar)$ has an empty interior, the only intervals included in $V(\Ebar \cap \Kcal_{\qrec})$ are isolated points, which shows that the limit $\lim_{j\to\infty} V(\vPpr(j)) $ exists. 

    We now prove that $\limsup_{j\to\infty} d(\vPpr(j), \Ebar\cap\Kcal_{\qrec}) = 0$. 
    Let $\Ncal \subset \Kcal_{\qrec}$ be an arbitrary neighborhood of $\Ebar \cap \Kcal_{\qrec}$. 
    From Lemma~\ref{lemma:2_4} there exist constants $\delta,\varepsilon, \lambda \in \R^+ $ such that for any $\delta'\in (0,\delta], \lambda' \in(0,\lambda]$, and $\eta>0$ one may find an integer $T$ and a sequence $\{\hvPpr(j)\}_{j\geq T}$ such that $\sup_{j \geq T} \norm{\vPpr(j) - \hvPpr(j)} \leq \delta',$
    and $\sup_{j\geq T} |V(\vPpr(j)) - V(\hvPpr(j))| \leq \eta$. Also, for any $j \geq T+1$
    \begin{align*}
         V(\hvPpr(j) \leq &V(\hvPpr(j-1)) - \step{j-1}\varepsilon 
         + (\eta + \step{j-1} \varepsilon) \indicator{V(\hvPpr(j-1))\in [V(\Ebar\cap\Kcal_{\qrec})]_{\alpha}}.         
    \end{align*}
     For $j \geq T$, define $\tau(j) := \inf\{k \geq 0 | \hvPpr(k+j) \in \Ncal\}$. 
    For any integer $p$, define ${\tau^{p}(j):= \min(\tau(j), p)}$. 
    Then
    \begin{align*}
        &V(\hvPpr(j+\tau^p(j))) - V(\hvPpr(j)) 
        = 
        \sum_{i=j+1}^{j+\tau^p(j)} \{V(\hvPpr(i)) - V(\hvPpr(i-1))\} 
         \leq -\varepsilon \sum_{i=j+1}^{j+\tau^p(j)} \step{i-1},
    \end{align*}
    with the convention that, for any sequence $\{a_i\}$ and any integer $l$, $\sum_{i=l+1}^l a_i = 0$. 
    Therefore, 
    \begin{align*}
    V(\vPpr(j+\tau^p(j))) - V(\vPpr(j))&= V(\vPpr(j+\tau^p(j))) - V(\hvPpr(j+\tau^p(j))) \cr 
    & +   V(\hvPpr(j+\tau^p(j)))- V(\hvPpr(j)) \cr
    & + V(\hvPpr(j))  -V(\vPpr(j))
    \leq 2 \eta -  \varepsilon \sum_{i=j+1}^{j+\tau^p(j)} \step{i-1}.
    \end{align*}
Since $\{V(\vPpr(j)\}$ converges, for any $\varepsilon' > 0 $ there exists ${ T' > T} $ such that, for all $j \geq T'$,
\begin{align*}
    -\varepsilon' < V(\vPpr(j+\tau^p(j))) - V(\vPpr(j)) 
    \leq 2 \eta - \varepsilon\!\!\!\! \sum_{i=j+1}^{j+\tau^p(j)}\!\!\!\! \step{i-1}.
\end{align*}

This implies that, for all $j \geq T'$ and all integer $p \geq 0$, 
    $\sum_{i=j+1}^{j+\tau^p(j)} \step{i-1} \leq C(\varepsilon',\eta):= \varepsilon^{-1} (\varepsilon'+2\eta)$.

 Since $\sum_{i=j+1}^{j+\tau(j)}\step{i-1} = \lim_{p\to\infty} \sum_{i=j+1}^{j+\tau^p(j)} \step{i-1}$ and $\sum_{i=1}^{\infty} \step{i-1} = \infty$, the previous relation implies that, for all $j \geq T'$, $\tau(j)<\infty$, and $\sum_{i=j+1}^{j+\tau(j)}\step{i-1}\leq C(\varepsilon',\eta)$. 

 For any integer $p$, we have
 \[
 \vPpr(j+p)-\vPpr(j) = \sum_{i=j+1}^{j+p} \step{i-1}\vf(\vPpr(i-1)) + \sum_{i=j+1}^{j+p} \step{i-1}\vdeterror(i),\]
 which implies that 
 \begin{align*}
     \norm{\vPpr(j+p) - \vPpr(j)} \leq \norm{\vf}_{\Kcal_{\qrec}} \sum_{i=j+1}^{j+p} \step{i-1} + \norm{\sum_{i=j+1}^{j+p} \step{i-1}\vdeterror(i)}.
 \end{align*}
 Applying this inequality for $j \geq T'$ and $p = \tau(j)$ and using that, by definition $\hvPpr(j+\tau(j))\in \Ncal$, 
\begin{align*}
    d(\vPpr(j),\Ncal) 
    &\leq 
    \norm{ \hvPpr(j+\tau(j)) - \vPpr(j+\tau(j))} +
    \norm{ \vPpr(j+\tau(j)) = \vPpr(j)} \cr 
    &\leq \delta' +\norm{\vf}_{\Kcal_{\qrec}} C(\varepsilon',\eta) + \norm{\sum_{i=j+1}^{j+\tau(j)} \step{i-1}\vdeterror(i)}.
\end{align*}
Since $\eta,\delta'$, and $\varepsilon'$ can be chosen arbitrarily small, and ${\limsup_{k\to\infty} \sup_{l\geq k} \norm{\sum_{i=k}^l \step{i-1}\vdeterror(i)} = 0}$, the latter inequality shows that $\lim_{j\to\infty} d(\vPpr(j),\Ncal) = 0$. 
Since $\Ncal$ is arbitrary, we  have 
$${\lim_{j\to\infty} d(\vPpr(j), \Ebar\cap\Kcal_{\qrec}) = 0}.$$
\end{pf*}

\end{document}